\documentclass[a4paper,reqno]{amsart}
\usepackage[all]{xy}           %For commutative diagrams
\usepackage{amssymb}           %For double arrows
\usepackage{hyperref}
\numberwithin{equation}{section}

\newtheorem{definition}{Definition}[section]

\newtheorem{theorem}[definition]{Theorem}
\newtheorem{proposition}[definition]{Proposition}
\newtheorem{corollary}[definition]{Corollary}
\newtheorem{remarkth}[definition]{Remark}

\newenvironment{remark}{\begin{remarkth}\upshape}{\hfill$\diamond$\end{remarkth}}

\renewcommand{\emph}[1]{{\bfseries\itshape{#1}}}

\newcommand{\R}{\mathbb{R}}      %Numeros reales
\newcommand{\F}{\mathbb{F}}

\newcommand{\lcf}{\lbrack\! \lbrack}
\newcommand{\rcf}{\rbrack\! \rbrack}
\newcommand\map[3]{#1\ \colon\ #2\longrightarrow#3}
\newcommand{\pai}[2]{\langle#1,#2\rangle}

\newcommand{\lvec}[1]{\overleftarrow{#1}}
\newcommand{\rvec}[1]{\overrightarrow{#1}}

\renewcommand{\d}{\mathrm{d}^\circ}
\newcommand{\qquand}{\qquad\text{and}\qquad}
\newcommand{\quand}{\quad\text{and}\quad}
\newcommand{\cinfty}[1]{C^\infty(#1)}

\newcommand{\I}{I\mkern-7muI}
\newcommand{\tr}{\operatorname{Tr}}
\newcommand{\e}{\mathrm{e}}

\begin{document}

\title{Discrete Lagrangian and Hamiltonian Mechanics on Lie
groupoids}

\author[J.\ C.\ Marrero]{Juan C.\ Marrero}
\address{Juan C.\ Marrero:
Departamento de Matem\'atica Fundamental, Facultad de
Ma\-te\-m\'a\-ti\-cas, Universidad de la Laguna, La Laguna,
Tenerife, Canary Islands, Spain} \email{jcmarrer@ull.es}

\author[D.\ Mart\'{\i}n de Diego]{David Mart\'{\i}n de Diego}
\address{D.\ Mart\'{\i}n de Diego:
Instituto de Matem\'aticas y F{\'\i}sica Fundamental, Consejo
Superior de Investigaciones Cient\'{\i}ficas, Serrano 123, 28006
Madrid, Spain} \email{d.martin@imaff.cfmac.csic.es}

\author[E.\ Mart\'{\i}nez]{Eduardo Mart\'{\i}nez}
\address{Eduardo Mart\'{\i}nez:
Departamento de Matem\'atica Aplicada, Facultad de Ciencias,
Universidad de Zaragoza, 50009 Zaragoza, Spain}
\email{emf@unizar.es}

\thanks{This work has been partially supported by MICYT (Spain)
Grants BMF 2003-01319, MTM 2004-7832 and BMF 2003-02532. D. Mart\'{\i}n
de  Diego and E. Mart\'{\i}nez acknowledge to Department of Fundamental
Mathematics (University of La Laguna) for the hospitality offered
to their during the period when this paper was prepared.}



\begin{abstract}
The purpose of this paper is to describe geometrically discrete
Lagrangian and Hamiltonian Mechanics on Lie groupoids. From a
variational principle we derive the discrete Euler-Lagrange
equations and we introduce a symplectic 2-section, which is
preserved by the Lagrange evolution operator. In terms of the
discrete Legendre transformations we define the Hamiltonian
evolution operator which is a symplectic map with respect to the
canonical symplectic 2-section on the prolongation of the dual of
the Lie algebroid of the given groupoid. The equations we get
include as particular cases the classical discrete Euler-Lagrange
equations, the discrete Euler-Poincar\'e and discrete
Lagrange-Poincar\'e equations. Our results can be important for
the construction of geometric integrators for continuous
Lagrangian systems.
\end{abstract}

\maketitle

\vspace{-20pt}

 \tableofcontents

\section{Introduction}
During the last decade, much effort has been devoted to
construction of geometric integrators for Lagrangian systems using
a discrete variational principle (see \cite{mawest} and references
therein). In particular, this effort has been concentrated for the
case of discrete Lagrangian functions $L$ on the cartesian product
$Q\times Q$ of a differentiable manifold. This cartesian product
plays the role of a ``discretized version" of the standard
velocity phase space $TQ$.  Applying a natural discrete
variational principle, one
 obtains a second order recursion operator $\xi: Q\times Q\longrightarrow Q\times Q$
assigning to each input pair  $(x,y)$ the output pair $(y,z)$.
When the discrete Lagrangian is an approximation of a continuous
Lagrangian function (more appropriately, when the discrete
Lagrangian approximates the integral action for $L$) we obtain a
numerical integrator which inherits some of the geometric
properties of the continuous Lagrangian (symplecticity, momentum
preservation). Although this type of geometric  integrators have
been mainly considered  for conservative systems, the extension to
geometric integrators for more involved situations  is relatively
easy, since, in some sense, many of the constructions  mimic the
corresponding ones for the continuous counterpart.  In this sense,
it  has been recently shown how  discrete variational mechanics
can include forced or dissipative systems, holonomic constraints,
explicitely time-dependent systems, frictional contact,
nonholonomic constraints, multisymplectic fields theories... All
these geometric integrators have demonstrated, in worked examples,
an exceptionally good longtime behavior and obviously this
research is of great interest for numerical and geometric
considerations (see \cite{Hair}).

On the other hand, Moser and Veselov~\cite{Mose} consider also
discrete Lagrangian systems evolving on a Lie group. All this
examples leads to A.\ Weinstein~\cite{weinstein} to study discrete
mechanics on Lie groupoids, which is a structure that includes as
particular examples the case of cartesian products  $Q\times Q$ as
well as Lie groups.

A Lie groupoid $G$ is a natural generalization of the concept of a
Lie group, where now not all  elements are composable. The product
$g_1 g_2$ of two elements is only defined on the set of composable
pairs $G_2=\{(g, h)\in G\times G\;|\; \beta(g)=\alpha(h)\}$ where
$\alpha: G\longrightarrow M$ and $\beta: G\longrightarrow M$ are
the source and target maps over a base manifold $M$. This concept
was introduced in differential geometry by Ch.\ Ereshmann in the
1950's. The infinitesimal version of a Lie groupoid $G$ is the Lie
algebroid $AG\longrightarrow M$, which is the restriction of the
vertical bundle of $\alpha$ to the submanifold of the identities.

We may thought a Lie algebroid $A$ over a manifold $M$, with
projection $\tau: A\to M$,  as a generalized version of the
tangent bundle to $M$. The geometry and dynamics on Lie algebroids
have  been extensively studied during the past years. In
particular, one of the authors of this paper (see~\cite{mart})
developed a geometric formalism of mechanics on Lie algebroids
similar to Klein's formalism~\cite{Klein} of the ordinary
Lagrangian mechanics and more recently a description of the
Hamiltonian dynamics on a Lie algebroid was given
in~\cite{LMM,Medina} (see also~\cite{PP}).

The key concept in this theory is the prolongation, ${\mathcal
P}^{\tau} A$, of the Lie algebroid over the fibred projection
$\tau$ (for the Lagrangian formalism) and the prolongation,
${\mathcal P}^{\tau^*} A$, over the dual fibred projection
$\tau^*: A^*\longrightarrow M$ (for the Hamiltonian formalism).
See \cite{LMM} for more details. Of course, when the Lie algebroid
is $A=TQ$ we obtain that ${\mathcal P}^{\tau} A=T(TQ)$ and
${\mathcal P}^{\tau^*} A=T(T^*Q)$, recovering the classical case.
An alternative approach, using the linear Poisson structure on
$A^*$ and the canonical isomorphism between $T^*A$ and $T^*A^*$
was discussed in \cite{GGU}.

Taking as starting point the  results by A.\ Weinstein
\cite{weinstein}, we elucidate in this paper the geometry of
Lagrangian systems on Lie groupoids and its Hamiltonian
counterpart.  Weinstein gave a variational derivation of the
discrete Euler-Lagrange equations for a Lagrangian $L: G\to \R$ on
a Lie groupoid $G$. We show that the appropriate space to develop
a geometric formalism for these equations is the Lie algebroid
${\mathcal P}^{\tau}G \equiv V\beta\oplus_G V\alpha\to G$ (see
section~3 for the definition of the Lie algebroid structure). Note
that $\mathcal{P}^\tau G$ is the total space of the prolongation
of the Lie groupoid $G$ over the vector bundle projection
$\tau:AG\to M$, and that the Lie algebroid of $\mathcal{P}^\tau G$
is just the prolongation $\mathcal{P}^\tau (AG)$ (the space were
the continuous Lagrangian Mechanics is developed). Using the Lie
algebroid structure of $\mathcal{P}^\tau G$ we may describe
discrete Mechanics on the Lie groupoid $G$. In particular,
\begin{itemize}
\item[--] We give a variational derivation of the discrete Euler-Lagrange equations:
\[
\lvec{X}({g})(L)-\rvec{X}({h})(L)=0
\]
for every section $X$ of $AG$, where the right or left arrow
denotes the induced right and left-invariant vector field on $G$.

\item[--] We introduce two Poincar\'e-Cartan 1-sections $\Theta^+_L$ and $\Theta^-_L$,
and an unique  Poincar\'e-Cartan 2-section, $\Omega_L$, on the Lie
algebroid $P^{\tau} G\longrightarrow G$.

\item[--] We study the discrete Lagrangian evolution operator $\xi: G\longrightarrow G$ and  its preservation properties. In particular, we prove that $({\mathcal P}^{\tau}\xi,\xi)^*\Omega_L=\Omega_L$, where ${\mathcal P}^{\tau}\xi$ is the natural prolongation of $\xi$ to ${\mathcal P}^{\tau}G$.

\item[--] Reduction theory is stablished in terms of morphisms of Lie groupoids.

\item[--] The associated Hamiltonian formalism is developed using the discrete Legendre transformations $\F^+L:G\to A^*G$ and $\F^-L:G\to A^*G$.

\item[--] A complete  characterization of the regularity of a Lagrangian on a Lie groupoid is given in terms of the symplecticity of $\Omega_L$ or, alternatively, in terms of the regularity of the discrete Legendre transformations. In particular, Theorem \ref{5.13} solves the question  posed by Weinstein~\cite{weinstein} about  the regularity  conditions for a discrete Lagrangian function on more general Lie groupoids than the cartesian product $Q\times Q$. In the regular case, we define the Hamiltonian evolution operator and we prove that it defines a symplectic map.

\item[--] We prove a Noether's theorem for discrete Mechanics on Lie groupoids.

\item[--] Finally, some illustrative examples are shown, for instance, discrete Mechanics on the cartesian product $Q\times Q$, on Lie groups (discrete Lie-Poisson equations), on action Lie groupoids (discrete Euler-Poincar\'e equations) and on gauge or Atiyah Lie groupoids (discrete Lagrange-Poincar\'e
equations).

\end{itemize}
We expect that the results of this paper could be relevant in the
construction of new geometric integrators, in particular, for the
numerical integration of dynamical systems with symmetry.

The paper is structured as follows. In Section 2 we review some
basic results on Lie algebroids and Lie groupoids. Section 3 is
devoted to study the Lie algebroid structure of the vector bundle
${\mathcal P}^{\tau}G\equiv V\beta\oplus_G V\alpha\ \to G$. The
main results of the paper appear in Section 4, where the geometric
structure of discrete Mechanics on Lie groupoids is given.
Finally, in Section 5, we study several examples of the theory.

\section{Lie algebroids and Lie groupoids}

\subsection{Lie algebroids}\label{section2.1}

A \emph{Lie algebroid} $A$ over a manifold $M$ is a real vector
bundle $\tau: A \to M$ together with a Lie bracket $\lcf\cdot,
\cdot\rcf$ on the space $\Gamma(\tau)$ of the global cross
sections of $\tau: A \to M$ and a bundle map $\rho:A \to TM$,
called \emph{the anchor map}, such that if we also denote by
$\rho:\Gamma(\tau)\to {\frak X}(M)$ the homomorphism of
$C^\infty(M)$-modules induced by the anchor map then
\begin{equation}\label{Leib}
 \lcf X,fY\rcf=f\lcf X,Y\rcf + \rho(X)(f)Y,
\end{equation}
 for $X,Y\in \Gamma(\tau)$ and $f\in C^\infty(M)$ (see \cite{Ma}).

If $X, Y, Z \in \Gamma(\tau)$ and $f \in C^{\infty}(M)$ then,
using (\ref{Leib}) and the fact that $\lcf\cdot, \cdot\rcf$ is a
Lie bracket, we obtain that
\begin{equation}\label{f1}
\lcf \lcf X, Y\rcf, fZ\rcf = f(\lcf X, \lcf Y, Z\rcf\rcf - \lcf Y,
\lcf X, Z\rcf\rcf) + [\rho( X), \rho(Y)](f) Z.
\end{equation}
On the other hand, from (\ref{Leib}), it follows that
\begin{equation}\label{f2}
\lcf \lcf X, Y\rcf, fZ\rcf = f \lcf \lcf X, Y\rcf, Z\rcf + \rho
\lcf X, Y \rcf (f) Z.
\end{equation}
Thus, using (\ref{f1}), (\ref{f2}) and the fact that $\lcf\cdot,
\cdot\rcf$ is a Lie bracket, we conclude that
\[
\rho \lcf X, Y \rcf = [\rho( X), \rho(Y)],
\]
that is, $\rho:\Gamma(\tau)\to {\frak X}(M)$ is a homomorphism
between the Lie algebras $(\Gamma(\tau),\lcf\cdot,\cdot\rcf)$ and
$({\frak X}(M),[\cdot,\cdot])$.

If $(A, \lcf\cdot, \cdot\rcf, \rho)$ is a Lie algebroid over $M$,
one may define \emph{the differential of $A$}, $d:\Gamma(\wedge^k
\tau^*)\to \Gamma(\wedge^{k+1}\tau^*)$, as follows
\begin{equation}\label{dif}
\begin{aligned}
d \mu(X_0,\dots, X_k)&=\sum_{i=0}^{k}
(-1)^i\rho(X_i)(\mu(X_0,\dots,
\widehat{X_i},\dots, X_k)) \\
&+ \sum_{i<j}(-1)^{i+j}\mu(\lcf X_i,X_j\rcf,X_0,\dots,
\widehat{X_i},\dots,\widehat{X_j},\dots ,X_k),
\end{aligned}
\end{equation}
for $\mu\in \Gamma(\wedge^k \tau^*)$ and $X_0,\dots ,X_k\in
\Gamma(\tau).$ $d$ is a cohomology operator, that is, $d^2=0$. In
particular, if $f: M\longrightarrow \R$ is a real smooth function
then $df(X)=\rho(X)f,$ for $X\in \Gamma(\tau)$. We may also define
the Lie derivative with respect to a section $X$ of $A$ as the
operator ${\mathcal L}_X: \Gamma(\Lambda^k A^*)\longrightarrow
\Gamma(\Lambda^k A^*)$ given by ${\mathcal L}_X=i_X\circ d+d\circ
i_X$ (for more details, see \cite{Ma}).

Trivial examples of Lie algebroids are a real Lie algebra ${\frak
g}$ of finite dimension (in this case, the base space is a single
point) and the tangent bundle $TM$ of a manifold $M$. Other
examples of Lie algebroids are: i) the vertical bundle
$(\tau_{P})_{|V\pi}: V\pi \to P$ of a fibration $\pi: P \to M$
(and, in general, the tangent vectors to a foliation of finite
dimension on a manifold $P$); ii) \emph{the Atiyah algebroid
associated with a principal $G$-bundle }(see \cite{LMM,Ma}); iii)
\emph{the prolongation ${\mathcal P}^{\pi}A$ of a Lie algebroid
$A$ over a fibration $\pi: P \to M$ }(see \cite{HM,LMM}) and iv)
\emph{the action Lie algebroid }$A\ltimes f$ over a map $f: M' \to
M$ (see \cite{HM,LMM}).

 Now, let $(A,\lcf\cdot,\cdot \rcf,\rho)$ (resp.,
$(A',\lcf\cdot,\cdot\rcf', \rho')$) be a Lie algebroid over a
manifold $M$ (resp., $M'$) and suppose that $\Psi: A \to A'$ is a
vector bundle morphism over the map $\Psi_{0}: M \to M'$. Then,
the pair $(\Psi, \Psi_{0})$ is said to be a \emph{Lie algebroid
morphism} if
\begin{equation}\label{dd'}
d ((\Psi, \Psi_{0})^*\phi')= (\Psi, \Psi_{0})^*(d'\phi'), \;\;\;
\text{ for all }\phi'\in \Gamma(\wedge^k(A')^*) \text{ and for all
}k,
\end{equation}
where $d$ (resp., $d'$) is the differential of the Lie algebroid
$A$ (resp., $A'$) (see \cite{LMM}). In the particular case when $M
= M'$ and $\Psi_{0} = Id$ then (\ref{dd'}) holds if and only if
\[
\lcf \Psi \circ X, \Psi \circ Y \rcf' = \Psi \lcf X, Y \rcf,
\makebox[.3cm]{} \rho'(\Psi X) = \rho(X), \makebox[.3cm]{} \mbox{
for } X, Y \in \Gamma(\tau).
\]

\subsection{Lie groupoids}\label{section2.2}

In this Section, we will recall the definition of a Lie groupoid
and some generalities about them are explained (for more details,
see \cite{CaWe,Ma}).

A \emph{groupoid} over a set $M$ is a set $G$ together with the
following structural maps:
\begin{itemize}
\item A pair of maps $\alpha: G \to M$, the \emph{source}, and
$\beta: G \to M$, the \emph{target}. Thus, an element $g \in G$ is
thought as an arrow from $x= \alpha(g)$ to $y = \beta(g)$ in $M$
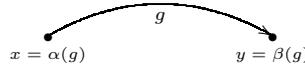
$$
\xymatrix{*=0{\stackrel{\bullet}{\mbox{\tiny
 $x=\alpha(g)$}}}{\ar@/^1pc/@<1ex>[rrr]_g}&&&*=0{\stackrel{\bullet}{\mbox{\tiny
$y=\beta(g)$}}}}
$$
The maps $\alpha$ and $\beta$ define the set of composable pairs
$$
G_{2}=\{(g,h) \in G \times G / \beta(g)=\alpha(h)\}.
$$
\item A \emph{multiplication} $m: G_{2} \to G$, to be denoted simply by $m(g,h)=gh$, such that
\begin{itemize}
\item $\alpha(gh)=\alpha(g)$ and $\beta(gh)=\beta(h)$.
\item $g(hk)=(gh)k$.
\end{itemize}
If $g$ is an arrow from $x = \alpha(g)$ to $y = \beta(g) =
\alpha(h)$ and $h$ is an arrow from $y$ to $z = \beta(h)$ then
$gh$ is the composite arrow from $x$ to $z$
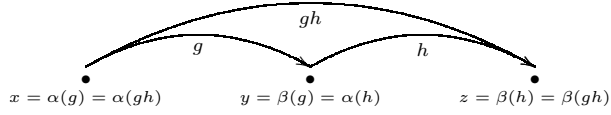
$$\xymatrix{*=0{\stackrel{\bullet}{\mbox{\tiny
 $x=\alpha(g)=\alpha(gh)$}}}{\ar@/^2pc/@<2ex>[rrrrrr]_{gh}}{\ar@/^1pc/@<2ex>[rrr]_g}&&&*=0{\stackrel{\bullet}{\mbox{\tiny
 $y=\beta(g)=\alpha(h)$}}}{\ar@/^1pc/@<2ex>[rrr]_h}&&&*=0{\stackrel{\bullet}{\mbox{\tiny
 $z=\beta(h)=\beta(gh)$}}}}$$
\item An \emph{identity section} $\epsilon: M \to G$ such that
\begin{itemize}
\item $\epsilon(\alpha(g))g=g$ and $g\epsilon(\beta(g))=g$.
\end{itemize}
\item An \emph{inversion map} $i: G \to G$, to be denoted simply by $i(g)=g^{-1}$, such that
\begin{itemize}
\item $g^{-1}g=\epsilon(\beta(g))$ and $gg^{-1}=\epsilon(\alpha(g))$.
\end{itemize}
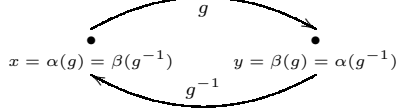
$$\xymatrix{*=0{\stackrel{\bullet}{\mbox{\tiny
 $x=\alpha(g)=\beta(g^{-1})$}}}{\ar@/^1pc/@<2ex>[rrr]_g}&&&*=0{\stackrel{\bullet}{\mbox{\tiny
 $y=\beta(g)=\alpha(g^{-1})$}}}{\ar@/^1pc/@<2ex>[lll]_{g^{-1}}}}$$

\end{itemize}

A groupoid $G$ over a set $M$ will be denoted simply by the symbol
$G \rightrightarrows M$.

The groupoid $G \rightrightarrows M$ is said to be a \emph{Lie
groupoid} if $G$ and $M$ are manifolds and all the structural maps
are differentiable with $\alpha$ and $\beta$ differentiable
submersions. If $G \rightrightarrows M$ is a Lie groupoid then $m$
is a submersion, $\epsilon$ is an immersion and $i$ is a
diffeomorphism. Moreover, if $x \in M$, $\alpha^{-1}(x)$ (resp.,
$\beta^{-1}(x)$) will be said the \emph{$\alpha$-fiber} (resp.,
the \emph{$\beta$-fiber}) of $x$.

On the other hand, if $g \in G$ then the \emph{left-translation by
$g \in G$} and the \emph{right-translation by $g$} are the
diffeomorphisms
$$
\begin{array}{lll}
l_{g}: \alpha^{-1}(\beta(g)) \longrightarrow
\alpha^{-1}(\alpha(g))&; \; \;& h \longrightarrow
l_{g}(h) = gh, \\
r_{g}: \beta^{-1}(\alpha(g)) \longrightarrow
\beta^{-1}(\beta(g))&; \; \;& h \longrightarrow r_{g}(h) = hg.
\end{array}
$$
Note that $l_{g}^{-1} = l_{g^{-1}}$ and $r_{g}^{-1} = r_{g^{-1}}$.

A vector field $\tilde{X}$ on $G$ is said to be
\emph{left-invariant} (resp., \emph{right-invariant}) if it is
tangent to the fibers of $\alpha$ (resp., $\beta$) and
$\tilde{X}(gh) = (T_{h}l_{g})(\tilde{X}_{h})$ (resp.,
$\tilde{X}(gh) = (T_{g}r_{h})(\tilde{X}(g)))$, for $(g,h) \in
G_{2}$.

Now, we will recall the definition of the \emph{Lie algebroid
associated with $G$}.

We consider the vector bundle $\tau: AG \to M$, whose fiber at a
point $x \in M$ is $A_{x}G = V_{\epsilon(x)}\alpha = Ker
(T_{\epsilon(x)}\alpha)$. It is easy to prove that there exists a
bijection between the space $\Gamma (\tau)$ and the set of
left-invariant (resp., right-invariant) vector fields on $G$. If
$X$ is a section of $\tau: AG \to M$, the corresponding
left-invariant (resp., right-invariant) vector field on $G$ will
be denoted $\lvec{X}$ (resp., $\rvec{X}$), where
\begin{equation}\label{linv}
\lvec{X}(g) = (T_{\epsilon(\beta(g))}l_{g})(X(\beta(g))),
\end{equation}
\begin{equation}\label{rinv}
\rvec{X}(g) = -(T_{\epsilon(\alpha(g))}r_{g})((T_{\epsilon
(\alpha(g))}i)( X(\alpha(g)))),
\end{equation}
for $g \in G$. Using the above facts, we may introduce a Lie
algebroid structure $(\lcf\cdot , \cdot\rcf, \rho)$ on $AG$, which
is defined by
\begin{equation}\label{LA}
\lvec{\lcf X, Y\rcf} = [\lvec{X}, \lvec{Y}], \makebox[.3cm]{}
\rho(X)(x) = (T_{\epsilon(x)}\beta)(X(x)),
\end{equation}
for $X, Y \in \Gamma(\tau)$ and $x \in M$. Note that
\begin{equation}\label{RL}
\rvec{\lcf X, Y\rcf} = -[\rvec{X}, \rvec{Y}], \makebox[.3cm]{}
[\rvec{X}, \lvec{Y}] = 0,
\end{equation}
\begin{equation}\label{2.6'}
Ti\circ \rvec{X}=-\lvec{X}\circ i,\;\;\;\; Ti\circ
\lvec{X}=-\rvec{X}\circ i,
\end{equation}
(for more details, see \cite{CDW,Ma}).

Given two Lie groupoids $G \rightrightarrows M$ and $G'
\rightrightarrows M'$, a \emph{morphism of Lie groupoids} is a
smooth map $\Phi: G \to G'$ such that
\[
(g, h) \in G_{2} \Longrightarrow (\Phi(g), \Phi(h)) \in (G')_{2}
\]
and
\[
\Phi(gh) = \Phi(g)\Phi(h).
\]
A morphism of Lie groupoids $\Phi: G \to G'$ induces a smooth map
$\Phi_{0}: M \to M'$ in such a way that
\[
\alpha' \circ \Phi = \Phi_{0} \circ \alpha, \makebox[.3cm]{}
\beta' \circ \Phi = \Phi_{0} \circ \beta, \makebox[.3cm]{} \Phi
\circ \epsilon = \epsilon' \circ \Phi_{0},
\]
$\alpha$, $\beta$ and $\epsilon$ (resp., $\alpha'$, $\beta'$ and
$\epsilon'$) being the source, the target and the identity section
of $G$ (resp., $G'$).

Suppose that $(\Phi, \Phi_{0})$ is a morphism between the Lie
groupoids $G \rightrightarrows M$ and $G' \rightrightarrows M'$
and that $\tau: AG \to M$ (resp., $\tau': AG' \to M'$) is the Lie
algebroid of $G$ (resp., $G'$). Then, if $x \in M$ we may consider
the linear map $A_{x}(\Phi): A_{x}G \to A_{\Phi_{0}(x)}G'$ defined
by
\begin{equation}\label{Amor}
A_{x}(\Phi)(v_{\epsilon(x)}) =
(T_{\epsilon(x)}\Phi)(v_{\epsilon(x)}), \; \; \mbox{ for }
v_{\epsilon(x)} \in A_{x}G.
\end{equation}
In fact, we have that the pair $(A(\Phi), \Phi_{0})$ is a morphism
between the Lie algebroids $\tau: AG \to M$ and $\tau': AG' \to
M'$ (see \cite{Ma}).

Next, we will present some examples of Lie groupoids.

1.- {\bf Lie groups}. Any Lie group $G$ is a Lie groupoid over
$\{\frak e \}$, the identity element of $G$. The Lie algebroid
associated with $G$ is just the Lie algebra ${\frak g}$ of $G$.

2.- {\bf The pair or banal groupoid}. Let $M$ be a manifold. The
product manifold $M \times M$ is a Lie groupoid over $M$ in the
following way: $\alpha$ is the projection onto the first factor
and $\beta$ is the projection onto the second factor; $\epsilon(x)
= (x, x)$, for all $x \in M$, $m((x, y), (y, z)) = (x, z)$, for
$(x, y), (y, z) \in M \times M$ and $i(x, y) = (y, x)$. $M \times
M \rightrightarrows M$ is called the \emph{pair or banal
groupoid}. If $x$ is a point of $M$, it follows that
\[
V_{\epsilon(x)}\alpha=\{0_x\}\times T_xM\subseteq T_xM\times T_xM
\cong T_{(x,x)}(M\times M).
\]
Thus, the linear maps
\[
\Psi_x:T_xM\to V_{\epsilon(x)}\alpha,\;\;\; v_x\to (0_x,v_x),
\]
induce an isomorphism (over the identity of $M$) between  the Lie
algebroids $\tau_M:TM\to M$ and $\tau:A(M\times M)\to M.$

3.- {\bf The Lie groupoid associated with a fibration}. Let $\pi:
P \to M$ be a fibration, that is, $\pi$ is a surjective submersion
and denote by $G_{\pi}$ the subset of $P \times P$ given by
\[
G_{\pi} = \{(p, p') \in P \times P / \pi(p) = \pi(p') \}.
\]
Then, $G_{\pi}$ is a Lie groupoid over $P$ and the structural maps
$\alpha_{\pi}$, $\beta_{\pi}$, $m_{\pi}$, $\epsilon_{\pi}$ and
$i_{\pi}$ are the restrictions to $G_{\pi}$ of the structural maps
of the pair groupoid $P \times P \rightrightarrows P$.

If $p$ is a point of $P$ it follows that
\[
V_{\epsilon_{\pi}(p)}\alpha_{\pi} = \{(0_{p}, Y_{p}) \in T_{p}P
\times T_{p}P / (T_{p}\pi)(Y_{p}) = 0 \}.
\]
Thus, if $(\tau_{P})_{|V\pi}: V\pi \to P$ is the vertical bundle
to $\pi$ then the linear maps
\[
(\Psi_{\pi})_{p}: V_{p}\pi \longrightarrow V_{\epsilon_{\pi}(p)}
\alpha_{\pi}, \makebox[.3cm]{} Y_{p} \longrightarrow (0_{p},
Y_{p})
\]
induce an isomorphism (over the identity of $M$) between the Lie
algebroids $(\tau_{P})_{|V\pi}: V\pi \to P$ and $\tau: AG_{\pi}
\to P$.

4.- {\bf Atiyah or gauge groupoids}. Let $p: Q \to M$ be a
principal $G$-bundle. Then, the free action, $\Phi: G \times Q \to
Q$, $(g, q) \to \Phi(g, q) = gq$, of $G$ on $Q$ induces, in a
natural way, a free action $\Phi \times \Phi: G \times (Q \times
Q) \to Q \times Q$ of $G$ on $Q \times Q$ given by $(\Phi \times
\Phi)(g, (q, q')) = (gq, gq')$, for $g \in G$ and $(q, q') \in Q
\times Q$. Moreover, one may consider the quotient manifold $(Q
\times Q) / G$ and it admits a Lie groupoid structure over $M$
with structural maps given by
\[
\begin{array}{lcl}
\tilde{\alpha}: (Q \times Q) / G \longrightarrow M &; \; \; &[(q,
q')] \longrightarrow p(q),
\\
\tilde{\beta}: (Q \times Q) / G \longrightarrow M &; \; \; &[(q,
q')] \longrightarrow p(q'),
\\
\tilde{\epsilon}: M \longrightarrow (Q \times Q) / G &; \; \;&
x\longrightarrow [(q, q)], \; \; \mbox{ if } p(q) = x,
\\
\tilde{m}: ((Q \times Q) / G)_{2} \longrightarrow (Q \times Q) / G
&; \; \;& ([(q, q')], [(gq', q'')]) \longrightarrow [(gq, q'')],
\\
\tilde{i}: (Q \times Q) / G \longrightarrow (Q \times Q) / G &; \;
\;& [(q, q')] \longrightarrow [(q', q)].
\end{array}
\]
This Lie groupoid is called \emph{the Atiyah (gauge) groupoid
associated with the principal $G$-bundle $p: Q \to M$} (see
\cite{L}).

If $x$ is a point of $M$ such that $p(q) = x$, with $q \in Q$, and
$p_{Q \times Q}: Q \times Q \to (Q \times Q) / G$ is the canonical
projection then it is clear that
\[
V_{\tilde{\epsilon}(x)}\tilde{\alpha} = (T_{(q, q)}p_{Q \times
Q})(\{0_{q}\} \times T_{q}Q).
\]
Thus, if $\tau_{Q}|G: TQ / G \to M$ is the Atiyah algebroid
associated with the principal $G$-bundle $p: G \to M$ then the
linear maps
\[
(TQ / G)_{x} \to V_{\tilde{\epsilon}(x)}\tilde{\alpha} \; \; ; \;
\; [v_{q}] \to (T_{(q, q)}p_{Q \times Q})(0_{q}, v_{q}),\mbox{
with $v_{q} \in T_{q}Q$},
\]
induce an isomorphism (over the identity of $M$) between the Lie
algebroids $\tau: A((Q \times Q) / G) \to M$ and $\tau_{Q} | G: TQ
/ G \to M$.

5.- {\bf The prolongation of a Lie groupoid over a fibration}.
Given a Lie groupoid $G \rightrightarrows M$ and a fibration $\pi:
P \to M$, we consider the set
\[
{\mathcal P}^{\pi}G = P \mbox{$\;$}_\pi \kern-3pt\times_\alpha G
\mbox{$\;$}_\beta \kern-3pt\times_\pi P = \{ (p, g, p') \in P
\times G \times P / \pi(p) = \alpha(g), \; \; \beta(g) = \pi(p')
\}.
\]
Then, ${\mathcal P}^{\pi}G$ is a Lie groupoid over $P$ with
structural maps given by
\[
\begin{array}{lcl}
\alpha^{\pi}: {\mathcal P}^{\pi}G \longrightarrow P &; \; \; &(p,
g, p') \longrightarrow p,
\\
\beta^{\pi}: {\mathcal P}^{\pi}G \longrightarrow P &; \; \;& (p,
g, p') \longrightarrow p',
\\
\epsilon^{\pi}: P \longrightarrow {\mathcal P}^{\pi}G &; \; \;&
p\longrightarrow (p, \epsilon(\pi(p)), p),
\\
m^{\pi}: ({\mathcal P}^{\pi}G)_{2} \longrightarrow {\mathcal
P}^{\pi}G &; \; \;& ((p, g, p'), (p', h, p'')) \longrightarrow (p,
gh, p''),
\\
i^{\pi}: {\mathcal P}^{\pi}G \longrightarrow {\mathcal P}^{\pi}G
&; \; \;& (p, g, p') \longrightarrow (p', g^{-1}, p).
\end{array}
\]
${\mathcal P}^{\pi}G$ is called the \emph{prolongation of $G$ over
$\pi: P \to M$}.

Now, denote by $\tau: AG \to M$ the Lie algebroid of $G$, by
$A({\mathcal P}^{\pi}G)$ the Lie algebroid of ${\mathcal P}^\pi G$
and by ${\mathcal P}^{\pi}(AG)$ the prolongation of $\tau: AG \to
M$ over the fibration $\pi$. If $p \in P$ and $m=\pi(p)$, then it
follows that
\[
A_{p}({\mathcal P}^{\pi}G) =  \{ (0_{p}, v_{\epsilon(m)}, X_{p})
\in T_{p}P \times A_{m}G \times T_{p}P / (T_{p}\pi)(X_{p}) =
(T_{\epsilon(m)}\beta)(v_{\epsilon(m)})\}
\]
%\[
%A_{p}({\mathcal P}^{\pi}G) =  \{ (0_{p}, v_{\epsilon(\pi(p))}, X_{p})
%\in T_{p}P \times A_{\pi(p)}G \times T_{p}P / (T_{p}\pi)(X_{p}) =
%(T_{\epsilon(\pi(p))}\beta)(v_{\epsilon(\pi(p))})\}
%\]
and, thus, one may consider the linear isomorphism
\begin{equation}\label{Isom}
(\Psi^{\pi})_{p}: A_{p}({\mathcal P}^{\pi}G) \longrightarrow
{\mathcal P}^{\pi}_{p}(AG), \makebox[.3cm]{} (0_{p},
v_{\epsilon(m)}, X_{p}) \longrightarrow (v_{\epsilon(m)}, X_{p}).
\end{equation}
%\begin{equation}\label{Isom}
%(\Psi^{\pi})_{p}: A_{p}({\mathcal P}^{\pi}G) \longrightarrow
%{\mathcal P}^{\pi}_{p}(AG), \makebox[.3cm]{} (0_{p},
%v_{\epsilon(\pi(p))}, X_{p}) \longrightarrow
%(v_{\epsilon(\pi(p))}, X_{p}).
%\end{equation}
In addition, one may prove that the maps $(\Psi^{\pi})_{p}$, $p
\in P$, induce an isomorphism $\Psi^{\pi}: A({\mathcal P}^{\pi}G)
\to {\mathcal P}^{\pi}(AG)$ between the Lie algebroids
$A({\mathcal P}^{\pi}G)$ and ${\mathcal P}^{\pi}(AG)$ (for more
details, see \cite{HM}).

6.- {\bf Action Lie groupoids}. Let $G \rightrightarrows M$ be a
Lie groupoid and $\pi: P \to M$ be a smooth map. If $P
\mbox{$\;$}_\pi \kern-3pt\times_\alpha G = \{ (p, g) \in P \times
G / \pi(p) = \alpha (g) \}$ then a right action of $G$ on $\pi$ is
a smooth map
\[
P \mbox{$\;$}_\pi \kern-3pt\times_\alpha G \to P, \; \; (p, g) \to
pg,
\]
which satisfies the following relations
$$\begin{array}{rcll}
\pi(pg) &=& \beta(g), &\mbox{ for } (p, g) \in P \mbox{$\;$}_\pi
\kern-3pt\times_\alpha G,\\
(pg)h &=& p(gh),    &\mbox{ for } (p, g) \in P \mbox{$\;$}_\pi
\kern-3pt\times_\alpha G  \mbox{ and } (g, h) \in G_{2}, \mbox{ and }\\
p\epsilon(\pi(p)) &=& p,   &\mbox{ for } p \in P.
\end{array}$$

Given such an action one constructs \emph{the action Lie groupoid}
$P \mbox{$\;$}_\pi \kern-3pt\times_\alpha G$ over $P$ by defining
\[
\begin{array}{lcl}
\tilde\alpha_{\pi}:  P \mbox{$\;$}_\pi \kern-3pt\times_\alpha G
\longrightarrow P &; \; \; &(p, g) \longrightarrow p,
\\
\tilde\beta_{\pi}: P \mbox{$\;$}_\pi \kern-3pt\times_\alpha G
\longrightarrow P &; \; \;& (p, g) \longrightarrow pg,
\\
\tilde\epsilon_{\pi}: P \longrightarrow P \mbox{$\;$}_\pi
\kern-3pt\times_\alpha G  &; \; \;& p \longrightarrow (p,
\epsilon(\pi(p))),
\\
\tilde{m}_{\pi}: (P \mbox{$\;$}_\pi \kern-3pt\times_\alpha G)_{2}
\longrightarrow P \mbox{$\;$}_\pi \kern-3pt\times_\alpha G &; \;
\;& ((p, g), (pg, h)) \longrightarrow (p, gh),
\\
\tilde{i}_{\pi}:  P \mbox{$\;$}_\pi \kern-3pt\times_\alpha G
\longrightarrow P \mbox{$\;$}_\pi \kern-3pt\times_\alpha G &; \;
\;& (p, g) \longrightarrow (pg, g^{-1}).
\end{array}
\]
Now, if $p \in P$, we consider the map $p \;\cdot :
\alpha^{-1}(\pi(p)) \to P$ given by
\[
p \cdot (g) = pg.
\]
Then, if $\tau: AG \to M$ is the Lie algebroid of $G$, the
$\R$-linear map $\Phi: \Gamma(\tau) \to {\frak X}(P)$ defined by
\[
\Phi(X)(p) = (T_{\epsilon(\pi(p))}p \;\cdot)(X(\pi(p))), \; \;
\mbox{ for } X \in \Gamma(\tau) \mbox{ and } p \in P,
\]
induces an action of $AG$ on $\pi: P \to M$. In addition, the Lie
algebroid associated with the Lie groupoid $P \mbox{$\;$}_\pi
\kern-3pt\times_\alpha G \rightrightarrows P$ is the action Lie
algebroid $AG \ltimes \pi$ (for more details, see \cite{HM}).

%%%%%%%%%%%%%%%%%%%%%%%%%%%%%%%%%%%%%%%%%%%%%%%%%%%%%%%%

\section{Lie algebroid structure on the vector bundle $\pi^\tau: {\mathcal P}^{\tau}G \to G$}

Let $G \rightrightarrows M$ be a Lie groupoid with structural maps
\[
\alpha, \beta: G \to M, \; \; \epsilon: M \to G, \; \; i: G \to G,
\; \; m: G_{2} \to G.
\]
Suppose that $\tau: AG \to M$ is the Lie algebroid of $G$ and that
${\mathcal P}^{\tau}G $ is the prolongation of $G$ over the
fibration $\tau: AG \to M$ (see Example $5$ in Section
\ref{section2.2}), that is,
\[
{\mathcal P}^{\tau}G = AG \mbox{$\;$}_\tau \kern-3pt\times_\alpha
G \mbox{$\;$}_\beta \kern-3pt\times_\tau AG.
\]
${\mathcal P}^{\tau}G$ is a Lie groupoid over $AG$ and we may
define the bijective map $\Theta: {\mathcal P}^{\tau}G \to V\beta
\oplus_{G} V\alpha$ as follows
\[
\Theta(u_{\epsilon(\alpha(g))}, g, v_{\epsilon(\beta(g))}) =
((T_{\epsilon(\alpha(g))}(r_{g} \circ
i))(u_{\epsilon(\alpha(g))}), (T_{\epsilon(\beta(g))}l_{g})
(v_{\epsilon(\beta(g))})),
\]
for $(u_{\epsilon(\alpha(g))}, g, v_{\epsilon(\beta(g))}) \in
A_{\alpha(g)}G \times G \times A_{\beta(g)}G$. Thus, $V\beta
\oplus_{G}V\alpha$ is a Lie groupoid over $AG$ (this Lie groupoid
was considered by Saunders \cite{Sau}). We remark that the Lie
algebroid of ${\mathcal P}^{\tau}G \equiv V\beta\oplus_{G}V\alpha
\rightrightarrows AG$ is isomorphic to the prolongation of $AG$
over $\tau: AG \to M$ and that the prolongation of a Lie algebroid
$A$ over the vector bundle projection $\tau: A \to M$ plays an
important role in the description of Lagrangian Mechanics on $A$
(see \cite{LMM,mart}).

On the other hand, note that ${\mathcal P}^{\tau}G \equiv V\beta
\oplus_{G}V\alpha$ is a real vector bundle over $G$. In this
section, we will prove that the vector bundle $\pi^\tau:{\mathcal
P}^{\tau}G \to G$ admits an integrable Lie algebroid structure. In
other words, we will prove that there exists a Lie groupoid $H
\rightrightarrows G$ over $G$ such that the Lie algebroid $AH$ is
isomorphic to the real vector bundle $\pi^\tau: {\mathcal
P}^{\tau}G \to G$. In addition, we will see that the Lie groupoid
$H$ is isomorphic to the prolongations of $G$ over $\alpha$ and
$\beta$.

It is clear that the Lie algebroids of the Lie groupoids over $G$
\[
G_{\beta} = \{(g, h) \in G \times G / \beta(g) = \beta(h)\}, \; \;
G_{\alpha} = \{(r, s) \in G \times G / \alpha(r) = \alpha(s)\},
\]
are just $V\beta \to G$ and $V \alpha \to G$, respectively. This
fact suggests to consider the following manifold
\[
G_{\beta} \star G_{\alpha} = \{ ((g, h), (r, s)) \in G_{\beta}
\times G_{\alpha} / \beta_{\beta} (g, h) = \alpha_{\alpha}(r,
s)\},
\]
where $\beta_{\beta}: G_{\beta} \to G$ (respectively,
$\alpha_{\alpha}: G_{\alpha} \to G$) is the target (respectively,
the source) of the Lie groupoid $G_{\beta} \rightrightarrows G$
(respectively, $G_{\alpha} \rightrightarrows G$).

We will identify the space $G_{\beta} \star G_{\alpha}$ with
\[
\{(g, h, s) \in G \times G \times G / \beta(g) = \beta(h), \;
\alpha(h) = \alpha(s)\}.
\]
This space admits a Lie groupoid structure over $G$ with
structural maps given by
\begin{equation}\label{e1}
\begin{array}{lll}
\alpha_{\beta \alpha}: G_{\beta} \star G_{\alpha} \longrightarrow
G &; \;
\;& \kern-12pt(g, h, s)  \longrightarrow g,\\
\beta_{\beta \alpha}: G_{\beta} \star G_{\alpha}  \longrightarrow
G &; \;
\;& \kern-12pt(g, h, s)  \longrightarrow s,\\
\epsilon_{\beta \alpha}: G  \longrightarrow G_{\beta} \star
G_{\alpha}& ; \;
\;& \kern-12pt g  \longrightarrow (g, g, g),\\
m_{\beta \alpha}: (G_{\beta} \star G_{\alpha})_{2} \longrightarrow
G_{\beta} \star G_{\alpha}& ; \; \;& \kern-12pt((g, h, s), (s, h',
s'))
 \longrightarrow (g, h's^{-1}h, s'),\\
i_{\beta \alpha}: G_{\beta} \star G_{\alpha}  \longrightarrow
G_{\beta} \star G_{\alpha}& ; \; \;& \kern-12pt(g, h, s)
\longrightarrow (s, gh^{-1}s, g).
\end{array}
\end{equation}
Note that
$$\begin{array}{lll} j_{\beta}: G_{\beta}
\longrightarrow G_{\beta} \star G_{\alpha} &; \; \;& (g, h)
\longrightarrow
j_{\beta}(g, h) = (g, h, h),\\
j_{\alpha}: G_{\alpha}  \longrightarrow G_{\beta} \star
G_{\alpha}& ; \; \; &(h, s)  \longrightarrow j_{\alpha}(h, s) =
(h, h, s),
\end{array}
$$
are Lie groupoid morphisms and that the map
\[
m_{\beta \alpha}(j_{\beta}, j_{\alpha}):  G_{\beta} \star
G_{\alpha} \to G_{\beta} \star G_{\alpha} ; \; \; (g, h, s) \to
m_{\beta \alpha}(j_{\beta}(g, h), j_{\alpha}(h, s))
\]
is just the identity map. This implies that $(G_{\beta},
G_{\alpha})$ is a \emph{matched pair of Lie groupoids} in the
sense of Mackenzie \cite{M1} (see also \cite{Mo}).

Denote by $(\lcf\cdot , \cdot\rcf, \rho)$ the Lie algebroid
structure on $\tau: AG \to M$.

\begin{theorem} \label{t3.1}
Let $A(G_{\beta} \star G_{\alpha}) \to G$ be the Lie algebroid of
the Lie groupoid $G_{\beta} \star G_{\alpha} \rightrightarrows G$.
Then:
\begin{enumerate}
\item
The vector bundles $A(G_{\beta} \star G_{\alpha}) \to G$ and
$\pi^\tau:{\mathcal P}^{\tau}G \cong V\beta \oplus_{G} V\alpha \to
G$ are isomorphic. Thus, the vector bundle $\pi^\tau:{\mathcal
P}^{\tau}G \cong V\beta \oplus_{G} V\alpha \to G$ admits a Lie
algebroid structure.

\item
The anchor map $\rho^{{\mathcal P}^{\tau}G}$ of
$\pi^\tau:{\mathcal P}^{\tau}G \cong V\beta \oplus_{G} V\alpha \to
G$ is given by
\begin{equation} \label{e1'}
 \rho^{{\mathcal P}^{\tau}G}(X_{g}, Y_{g}) =
X_{g} + Y_{g}, \; \; \mbox{ for } (X_{g}, Y_{g}) \in V_{g}\beta
\oplus V_{g}\alpha,
\end{equation}
and the Lie bracket $\lcf\cdot , \cdot\rcf^{{\mathcal P}^{\tau}G}$
on the space $\Gamma(\pi^\tau)$ is characterized by the following
relation
\begin{equation}\label{e1''}
\lcf (\rvec{X}, \lvec{Y}), (\rvec{X'}, \lvec{Y'}) \rcf ^{{\mathcal
P}^{\tau}G} = (-\overrightarrow{\lcf X, X' \rcf},
\overleftarrow{\lcf Y, Y' \rcf}),
\end{equation}
for $X, Y, X', Y' \in \Gamma(\tau)$.
\end{enumerate}
\end{theorem}
\begin{proof}
(i) If $g \in G$ then, from (\ref{e1}), we deduce that the vector
space $V_{\epsilon_{\beta \alpha}(g)}\alpha_{\beta \alpha}$ may be
described as follows
\begin{align*}
V_{\epsilon_{\beta \alpha}(g)}\alpha_{\beta\alpha} &=
\{(0_{g},X_{g}, Z_{g}) \in T_{g}G \times T_{g}G \times T_{g}G /
X_{g} \in V_{g}\beta, \; (T_{g}\alpha)(X_{g}) =
(T_{g}\alpha)(Z_{g}) \}
\\
&\cong \{(X_{g}, Z_{g}) \in T_{g}G \times T_{g}G / X_{g} \in
V_{g}\beta, \; (T_{g}\alpha)(X_{g}) = (T_{g}\alpha)(Z_{g}) \}.
\end{align*}
Now, we will define the linear map $\Psi_{g}: V_{\epsilon_{\beta
\alpha}(g)}\alpha_{\beta\alpha} \to V_{g}\beta \oplus V_{g}\alpha
\cong {\mathcal P}^{\tau}_{g}G$ by
\begin{equation}\label{e2}
\Psi_{g}(X_{g}, Z_{g}) = (X_{g}, Z_{g} - X_{g}).
\end{equation}
It is clear that $\Psi_{g}$ is a linear isomorphism and
\begin{equation}\label{e3}
\Psi_{g}^{-1}(X_{g}, Y_{g}) = (X_{g}, X_{g} + Y_{g}), \; \; \mbox{
for } (X_{g}, Y_{g}) \in V_{g}\beta \oplus V_{g}\alpha \cong
{\mathcal P}^{\tau}_{g}G.
\end{equation}
Therefore, the collection of the maps $\Psi_{g}$, $g \in G$,
induces a vector bundle isomorphism $\Psi: A(G_{\beta} \star
G_{\alpha}) \to {\mathcal P}^{\tau}G \cong V\beta \oplus_{G}
V\alpha$ over the identity of $G$.

\medskip

(ii) A direct computation, using (\ref{e1}), proves that the
linear map $T_{\epsilon_{\beta \alpha}(g)}\beta_{\beta \alpha}:
V_{\epsilon_{\beta \alpha}(g)}\alpha_{\beta\alpha} \to T_{g}G$ is
given by
\begin{equation}\label{e4}
(T_{\epsilon_{\beta \alpha}(g)}\beta_{\beta \alpha})(X_{g}, Z_{g})
= Z_{g}.
\end{equation}
Consequently, from (\ref{LA}), (\ref{e3}) and (\ref{e4}), we
deduce that (\ref{e1'}) holds.

Next, we will prove (\ref{e1''}).

Using (\ref{e3}), it follows that
\[
(\Psi^{-1} \circ (\rvec{X}, \lvec{Y}))(g) = (0_{g}, \rvec{X}(g),
\rvec{X}(g) + \lvec{Y}(g)) \cong (\rvec{X}(g), \rvec{X}(g) +
\lvec{Y}(g)),
\]
for $g \in G$. Denote by $\lvec{\Psi^{-1}\circ (\rvec{X},
\lvec{Y})}$ the corresponding left-invariant vector field on
$G_{\beta} \star G_{\alpha}$. Then, from (\ref{linv}) and
(\ref{e1}), we have that
\[
\lvec{\Psi^{-1}\circ (\rvec{X}, \lvec{Y})}(g, h, s) = (0_{g},
\rvec{X}(h), \rvec{X}(s) + \lvec{Y}(s)), \; \; \mbox{ for } (g, h,
s) \in G_{\beta} \star G_{\alpha}.
\]
Thus, using (\ref{LA}) and (\ref{RL}), we conclude that
\[
[\lvec{\Psi^{-1}\circ (\rvec{X}, \lvec{Y})}, \lvec{\Psi^{-1}\circ
(\rvec{X'}, \lvec{Y'})}] = \lvec{\Psi^{-1}\circ(-\rvec{\lcf X,
X'\rcf}, \lvec{\lcf Y, Y'\rcf})}.
\]
Therefore, we obtain that (\ref{e1''}) holds.
\end{proof}

The above diagram shows the Lie groupoid and Lie algebroid
structures  of ${\mathcal P}^{\tau}G$:
\[
\xymatrix{%
&{\mathcal P}^{\tau}G\ar@{=>}[rr]^{\ \ \ \alpha^{\tau}}_{\ \ \
\beta^{\tau}}\ar[rd]^{\rho^{{\mathcal P}^{\tau}G}}
\ar[dd]^{\pi^\tau}&&AG\ar[rd]^{\rho}\ar'[d][dd]^{\tau}\\
&&TG\ar@{=>}[rr]^(.4){T\alpha}_(.4){T\beta}\ar[ld]_(.6){\tau_G}&&TM\ar[ld]^{\tau_M}\\
&G\ar@{=>}[rr]^\alpha_\beta&&M }
\]

Given a section $X$ of $AG\longrightarrow M$, we define the
sections $X^{(1, 0)}$, $X^{(0,1)}$ (the $\beta$ and $\alpha$-
lifts) and $X^{(1, 1)}$ (the complete lift)  of $X$ to $\pi^\tau:
{\mathcal P}^{\tau}G\longrightarrow G$ as follows:
\[
X^{(1, 0)}(g)=(\rvec{X}(g), 0_g), \quad X^{(0, 1)}(g)=(0_g,
\lvec{X}(g)) \quand X^{(1, 1)}(g)=(-\rvec{X}(g), \lvec{X}(g))
\]
We can easily see that
\begin{equation}\label{e1''*}
%\left\{
\begin{array}{l}
\lcf X^{(1, 0)}, Y^{(1, 0)}\rcf ^{{\mathcal P}^{\tau}G} = -\lcf X,
Y \rcf^{(1,0)} \\
 \lcf X^{(0, 1)}, Y^{(0, 1)}\rcf ^{{\mathcal P}^{\tau}G} =\hphantom{-} \lcf X, Y \rcf^{(0,1)}
\end{array}
%\right.
 \hbox{  and  } \ \lcf X^{(0,
1)}, Y^{(1, 0)}\rcf ^{{\mathcal P}^{\tau}G} = 0
\end{equation}
and, as a consequence,
\begin{equation}\label{e2''}
%\left\{
\begin{array}{l}
\lcf X^{(1, 1)}, Y^{(1, 0)}\rcf ^{{\mathcal P}^{\tau}G} = \lcf X,
Y \rcf^{(1,0)}  \\
\lcf X^{(1, 1)}, Y^{(0, 1)}\rcf ^{{\mathcal P}^{\tau}G} = \lcf X,
Y \rcf ^{(0,1)}
\end{array}
%\right.
 \hbox{  and  }\ \lcf X^{(1,
1)}, Y^{(1, 1)}\rcf ^{{\mathcal P}^{\tau}G} = \lcf X, Y
\rcf^{(1,1)}.
\end{equation}

\begin{remark}
From Theorem \ref{t3.1}, we deduce that the canonical inclusions
\[
(Id, 0): V\beta \to {\mathcal P}^{\tau}G \cong V\beta \oplus_{G}
V\alpha, \makebox[.3cm]{} (0, Id): V\alpha \to {\mathcal
P}^{\tau}G \cong V\beta \oplus_{G} V\alpha,
\]
are Lie algebroid morphisms over the identity of $G$. In other
words, $(V\beta, V\alpha)$ is a \emph{matched pair of Lie
algebroids} in the sense of Mokri \cite{Mo}. This fact directly
follows using the following general theorem (see \cite{Mo}): if
$(G, H)$ is a matched pair of Lie groupoids then $(AG, AH)$ is a
matched pair of Lie algebroids.
\end{remark}
Next, we will consider the prolongation ${\mathcal P}^{\beta}G$ of
the Lie groupoid $G$ over the target $\beta: G \to M$. We recall
that
\[
{\mathcal P}^{\beta}G = G \mbox{$\;$}_\beta \kern-3pt\times_\alpha
G \mbox{$\;$}_\beta \kern-3pt\times_\beta G = \{(g, h, s)\in G
\times G \times G / \beta(g) = \alpha(h), \; \beta(h) =
\beta(s)\},
\]
and that ${\mathcal P}^{\beta}G$ is a Lie groupoid over $G$ with
structural maps
\begin{equation}\label{e5}
\begin{array}{lll}
 \alpha^{\beta}: {\mathcal P}^{\beta}G
 \longrightarrow G &; \; \; &(g, h, s)  \longrightarrow g,\\
\beta^{\beta}: {\mathcal P}^{\beta}G  \longrightarrow G &; \; \;&
(g, h, s)  \longrightarrow
s,\\
\epsilon^{\beta}: G  \longrightarrow {\mathcal P}^{\beta}G&; \; \;
&g
 \longrightarrow (g,
\epsilon(\beta(g)), g),\\
m^{\beta}: ({\mathcal P}^{\beta}G)_{2}  \longrightarrow {\mathcal
P}^{\beta}G &; \; \;&
((g, h, s), (s, t, u))  \longrightarrow (g, ht, u),\\
i^{\beta}: {\mathcal P}^{\beta}G  \longrightarrow {\mathcal
P}^{\beta}G &; \; \;& (g, h, s)  \longrightarrow (s, h^{-1}, g).
\end{array}
\end{equation}

Moreover, we also have that the Lie algebroid of ${\mathcal
P}^{\beta}G$ may be identified with the prolongation ${\mathcal
P}^{\beta}(AG)$ of $AG$ over $\beta: G \to M$. We remark that
\[
{\mathcal P}^{\beta}_{g}(AG) = \{(v_{\epsilon(\beta(g))}, X_{g})
\in A_{\beta(g)}G \times T_{g}G /
(T_{\epsilon(\beta(g))}\beta)(v_{\epsilon(\beta(g))}) =
(T_{g}\beta)(X_{g}) \}
\]
for $g \in G$.

\begin{theorem}\label{t3.3}
Let $\Phi^{\beta}: G_{\beta} \star G_{\alpha} \to {\mathcal
P}^{\beta}G$ be the map defined by
\begin{equation}\label{e6}
\Phi^{\beta}(g, h, s) = (g, h^{-1}s, s),
\end{equation}
for $(g, h, s) \in G_{\beta} \star G_{\alpha}$. Then:
\begin{enumerate}
\item
$\Phi^{\beta}$ is a Lie groupoid isomorphism over the identity of
$G$.

\item
If $A(\Phi^{\beta}): A(G_{\beta} \star G_{\alpha}) \to A({\mathcal
P}^{\beta}G)$ is the corresponding Lie algebroid isomorphism then,
under the identifications
\[
A(G_{\beta} \star G_{\alpha}) \cong {\mathcal P}^{\tau}G \cong
V\beta \oplus_{G} V\alpha, \makebox[.3cm]{} A({\mathcal
P}^{\beta}G) \cong {\mathcal P}^{\beta}(AG),
\]
$A(\Phi^{\beta})$ is given by
\begin{equation} \label{e7}
A_g(\Phi^{\beta})(X_{g}, Y_{g}) = ((T_{g}l_{g^{-1}})(Y_{g}), X_{g}
+ Y_{g}),
\end{equation}
for $(X_{g}, Y_{g}) \in V_{g}\beta \oplus V_{g}\alpha$, where
$l_{g^{-1}}: \alpha^{-1}(\alpha(g)) \to \alpha^{-1}(\beta(g))$ is
the left-translation by $g^{-1}$.
\end{enumerate}
\end{theorem}
\begin{proof}
(i) A direct computation, using (\ref{e1}) and (\ref{e5}), proves
the result.

(ii) If $g \in G$ we have that
\[
\begin{array}{lcl} A_{g}(G_{\beta} \star G_{\alpha}) = V_{\epsilon_{\beta
\alpha}(g)}\alpha_{\beta \alpha} &=& \{(0_{g}, X_{g}, Z_{g}) \in
T_{g}G \times T_{g}G \times T_{g}G / X_{g} \in V_{g}\beta,\\ &&
(T_{g}\alpha)(X_{g}) = (T_{g}\alpha)(Z_{g}) \},
\\
A_{g}({\mathcal P}^{\beta}G) =
V_{\epsilon^{\beta}(g)}\alpha^{\beta} &=& \{(0_{g},
v_{\epsilon(\beta(g))}, Y_{g}) \in T_{g}G \times A_{\beta(g)}G
\times T_{g}G / \\&&(T_{g}\beta)(Y_{g}) =
(T_{\epsilon(\beta(g))}\beta)(v_{\epsilon(\beta(g))}) \}.
\end{array}
\]
Now, if $(0_{g}, X_{g}, Z_{g}) \in V_{\epsilon_{\beta
\alpha}(g)}\alpha_{\beta \alpha}$ then, from (\ref{e6}), we deduce
that
\begin{align*}
(T_{\epsilon_{\beta \alpha}(g)}\Phi^{\beta})&(0_{g}, X_{g}, Z_{g})
= (T_{\epsilon_{\beta \alpha}(g)}\Phi^{\beta})(0_{g}, 0_{g}, Z_{g}
- X_{g}) + (T_{\epsilon_{\beta \alpha}(g)}\Phi^{\beta})(0_{g},
X_{g}, X_{g})
\\
&= (0_{g}, (T_{g}l_{g^{-1}})(Z_{g} - X_{g}), Z_{g} - X_{g}) +
(T_{\epsilon_{\beta \alpha}(g)}\Phi^{\beta})(0_{g}, X_{g}, X_{g}).
\end{align*}
On the other hand, suppose that $\beta(g) = x \in M$ and that
$\gamma: (-\varepsilon, \varepsilon) \to \beta^{-1}(x)$ is a curve
in $\beta^{-1}(x)$ such that $\gamma(0) = g$ and $\gamma'(0) =
X_{g}$. Then, one may consider the curve $\tilde{\gamma}:
(-\varepsilon, \varepsilon) \to G_{\beta} \star G_{\alpha}$ on
$G_{\beta} \star G_{\alpha}$ given by
\[
\tilde{\gamma}(t) = (g, \gamma(t), \gamma(t))
\]
and it follows that
\[
\tilde{\gamma}(0) = (g, g, g), \makebox[.4cm]{} \tilde{\gamma}'(0)
= (0_{g}, X_{g}, X_{g}).
\]
Moreover, we obtain that
\[
\bar{\gamma}(t) = (\Phi^{\beta} \circ \tilde{\gamma})(t) = (g,
\epsilon(x), \gamma(t)), \makebox[.3cm]{} \mbox{ for all } t
\]
and thus
\[
\bar{\gamma}'(0) = (0_{g}, 0_{\epsilon(\beta(g))}, X_{g}).
\]
This proves that
\begin{equation}\label{e8}
(T_{\epsilon_{\beta \alpha}(g)}\Phi^{\beta})(0_{g}, X_{g}, Z_{g})
= (0_{g}, (T_{g}l_{g^{-1}})(Z_{g} - X_{g}), Z_{g}).
\end{equation}
Finally, using (\ref{Amor}), (\ref{Isom}), (\ref{e3}) and
(\ref{e8}), we deduce that (\ref{e7}) holds.
\end{proof}

Next, we will consider the prolongation ${\mathcal P}^{\alpha}G$
of the Lie groupoid $G$ over the source $\alpha: G \to M$. We
recall that
\[
{\mathcal P}^{\alpha}G = G \mbox{$\;$}_\alpha
\kern-3pt\times_\alpha G \mbox{$\;$}_\beta \kern-3pt\times_\alpha
G = \{(g, h, s) \in G \times G \times G / \alpha(g) = \alpha(h),
\beta(h) = \alpha(s) \}
\]
and that ${\mathcal P}^{\alpha}G$ is a Lie groupoid over $G$ with
structural maps
\[
\begin{array}{lcl}
\alpha^{\alpha}: {\mathcal P}^{\alpha}G \longrightarrow G &; \; \;
&(g, h, s) \longrightarrow g,
\\
\beta^{\alpha}: {\mathcal P}^{\alpha}G \longrightarrow G &; \; \;&
(g, h, s) \longrightarrow
s,\\
\epsilon^{\alpha}: G \longrightarrow {\mathcal P}^{\alpha}G &; \;
\;& g\longrightarrow (g,
\epsilon(\alpha(g)), g),\\
m^{\alpha}: ({\mathcal P}^{\alpha}G)_{2} \longrightarrow {\mathcal
P}^{\alpha}G &; \; \;& ((g, h, s),
(s, t, u)) \longrightarrow (g, ht, u),\\
i^{\alpha}: {\mathcal P}^{\alpha}G \longrightarrow {\mathcal
P}^{\alpha}G &; \; \;& (g, h, s) \longrightarrow (s, h^{-1}, g).
\end{array}
\]

Moreover, we also have that the Lie algebroid of ${\mathcal
P}^{\alpha}G$ may be identified with the prolongation ${\mathcal
P}^{\alpha}(AG)$ of $AG$ over $\alpha: G \to M$. We remark that
\[
{\mathcal P}^{\alpha}_{g}(AG) = \{ (v_{\epsilon(\alpha(g))},
X_{g}) \in A_{\alpha(g)}G \times T_{g}G / (T_{\epsilon(\alpha(g))}
\beta)(v_{\epsilon(\alpha(g))}) = (T_{g}\alpha)(X_{g}) \},
\]
for $g \in G$.
\begin{theorem}\label{t3.4}
Let $\Phi^{\alpha}: G_{\beta} \star G_{\alpha} \to {\mathcal
P}^{\alpha}G$ be the map defined by
\[
\Phi^{\alpha}(g, h, s) = (g, gh^{-1}, s),
\]
for $(g, h, s) \in G_{\beta} \star G_{\alpha}$. Then:
\begin{enumerate}
\item
$\Phi^{\alpha}$ is a Lie groupoid isomorphism over the identity of
$G$.

\item
If $A(\Phi^{\alpha}): A(G_{\beta} \star G_{\alpha}) \to
A({\mathcal P}^{\alpha}G)$ is the corresponding Lie algebroid
isomorphism then, under the canonical identifications
\[
A(G_{\beta} \star G_{\alpha}) \cong {\mathcal P}^{\tau}G \cong
V\beta \oplus_{G} V\alpha, \makebox[.4cm]{} A({\mathcal
P}^{\alpha}G) \cong {\mathcal P}^{\alpha}(AG),
\]
$A(\Phi^{\alpha})$ is given by
\begin{equation}\label{3.13}
A_g(\Phi^{\alpha})(X_{g}, Y_{g}) = (T_{g}(i \circ
r_{g^{-1}})(X_{g}), X_{g} + Y_{g}), \end{equation} for $(X_{g},
Y_{g}) \in V_{g}\beta \oplus V_{g}\alpha$, where $r_{g^{-1}}:
\beta^{-1}(\beta(g)) \to \beta^{-1}(\alpha(g))$ is the
right-translation by $g^{-1}$.

\end{enumerate}

\end{theorem}

\begin{proof}
Proceeding as in the proof of Theorem \ref{t3.3}, we deduce the
result.
\end{proof}

\section{Mechanics on  Lie Groupoids}

In this section, we introduce Lagrangian (Hamiltonian) Mechanics
on an arbitrary Lie groupoid and we will also analyze its
geometrical properties. This construction may be considered as a
discrete version of the construction of the Lagrangian
(Hamiltonian) Mechanics on Lie algebroids proposed in \cite{mart}
(see also \cite{LMM, Medina}). We first discuss discrete
Euler-Lagrange equations following a similar approach to
\cite{weinstein}, using a variational procedure. Secondly, we
intrinsically define and discuss the discrete Poincar\'e-Cartan
sections, Legendre transformations, regularity of the Lagrangian
and Noether's theorem.

\subsection{Discrete Euler-Lagrange equations}\label{section4.1}

Let $G$ be  a Lie groupoid with structural maps
\[
\alpha, \beta: G \to M, \; \; \epsilon: M \to G, \; \; i: G \to G,
\; \; m: G_{2} \to G.
\]
Denote by $\tau:AG\to M$ the Lie algebroid of $G$.

 A \emph{discrete Lagrangian} is a function $L: G \longrightarrow \R$. Fixed $g\in
G$, we define the set of admissible sequences with values in $G$:
\[
{\mathcal C}^N_{g}=\{(g_1, \ldots, g_N)\in G^N\; / \; (g_k,
g_{k+1})\in G_2 \hbox{ for } k=1,\ldots, N-1 \hbox{ and } g_1
\ldots g_n=g  \}.
\]

Given a  tangent vector  at $(g_1, \ldots, g_N)$ to the manifold
${\mathcal C}^N_{g}$, we may write it as the tangent vector at
$t=0$ of a curve in ${\mathcal C}^N_{g}$, $t\in (-\varepsilon,
\varepsilon)\subseteq \R\longrightarrow c(t)$ which  passes
through $(g_1, \ldots, g_N)$ at $t=0$. This type of curves is of
the form
\[
c(t)=(g_1h_1(t), h_1^{-1}(t)g_2 h_2(t), \ldots,
h_{N-2}^{-1}(t)g_{N-1} h_{N-1}(t), h_{N-1}^{-1}(t)g_N)
\]
where $h_k(t)\in \alpha^{-1}(\beta(g_k)),$ for all $t,$ and
$h_k(0)=\epsilon(\beta(g_k))$ for $k=1,\ldots, N-1$.

Therefore, we may  identify the   tangent space to ${\mathcal
C}^N_g$ at $(g_1,\dots ,g_N)$ with
\[
T_{(g_1, \ldots, g_N)}{\mathcal C}^N_g\equiv\left\{(v_1, \ldots,
v_{N-1})\; /\; v_k\in A_{x_k}G\hbox{ and } x_k=\beta(g_k), 1\leq
k\leq N-1\right\}.
\]
Observe that  each $v_k$ is the tangent vector to the
$\alpha$-vertical curve $h_k$ at $t=0$.

The curve $c$ is called a \emph{variation} of $(g_1, \ldots, g_N)$
and $(v_1, v_2, \ldots, v_{N-1})$ is called an \emph{infinitesimal
variation} of $(g_1, \ldots, g_N)$.

Define the  \emph{discrete action sum} associated to the discrete
Lagrangian $L: G\longrightarrow \R$
\[
\begin{array}{rcrcl}
{\mathcal S} L&:& {\mathcal C}^N_g& \longrightarrow &\R\\
          & &(g_1, \ldots, g_N)&\longmapsto&\displaystyle{ \sum_{k=1}^{N}
          L(g_k)}.
\end{array}
\]

We now proceed, as in the continuous case, to derive the discrete
equations of motion applying Hamilton's principle of critical
action. For it, we consider variations of the discrete action sum.

\begin{definition}[Discrete Hamilton's principle~\cite{weinstein}]
Given $g\in G$, an admissible sequence $(g_1, \ldots, g_N)\in
{\mathcal C}^N_g$ is a solution of the Lagrangian system
determined by $L:G\longrightarrow\R$ if and only if  $(g_1,
\ldots, g_N)$ is a critical point of ${\mathcal S}L$.
\end{definition}

Fist of all, in order to characterize the critical points, we need
to calculate: \begin{eqnarray*} \frac{d}{dt}\Big|_{t=0}{\mathcal
S} L(c(t))&=&\frac{d}{dt}\Big|_{t=0} \left\{
L(g_1h_1(t))+L(h_1^{-1}(t)g_2 h_2(t))\right. \\
&&\left. + \ldots + L(h_{N-2}^{-1}(t)g_{N-1} h_{N-1}(t))+
L(h^{-1}_{N-1}(t)g_N) \right\}.
 \end{eqnarray*}
 Therefore,
\begin{eqnarray*}
\frac{d}{dt}\Big|_{t=0}{\mathcal S} L(c(t))&=&
\sum_{k=1}^{N-1}\left(\d (L\circ l_{g_k})(\epsilon (x_k))(v_k)+\d
(L\circ r_{g_{k+1}}\circ i)(\epsilon(x_{k}))(v_k)\right)
\end{eqnarray*}
where $\d$ is the standard differential on $G$, i.e., $\d$ is the
differential of the Lie algebroid $\tau_G:TG\to G.$ Since the
critical condition is $\displaystyle{
\frac{d}{dt}\Big|_{t=0}{\mathcal S} L(c(t))=0}$ then, applying
(\ref{linv}) and (\ref{rinv}), we may rewrite this condition as
\[
0=\sum_{k=1}^{N-1}\left[\lvec{X}_k\big({g_k})(L)-\rvec{X}_k\big({g_{k+1}})(L)
\right] =\sum_{k=1}^{N-1}
\left[\pai{dL}{X^{(0,1)}_k}(g_k)-\pai{dL}{X^{(1,0)}_k}(g_{k+1})\right]
\]
where $d$ is the differential of the Lie algebroid
$\pi^\tau:{\mathcal P}^\tau G\equiv V\beta\oplus_G
V\alpha\longrightarrow G$ and $X_k$ is a section of $\tau:AG\to M$
such that $X_k(x_k)=v_k.$

For $N=2$ we obtain that $(g_1, g_2)\in G_2$ (with
$\beta(g_1)=\alpha(g_2)=x$) is a solution if
\[
\d\left[L\circ l_{g_1}+L\circ r_{g_2}\circ i\right](\epsilon
(x))_{| A_x G}=0
\]
or, alternatively,
\[
\lvec{X}({g_1})(L)-\rvec{X}({g_{2}})(L)=0
\]
for every section $X$ of $AG$. These equations will be called
\emph{discrete Euler-Lagrange equations}.

Thus, we may define the \emph{discrete Euler-Lagrange operator}:
\[
D_{\hbox{\footnotesize DEL}}L: G_2\longrightarrow A^*G\; ,
\]
where $A^*G$ is the dual of $AG$. This operator is given by
\[
D_{\hbox{\footnotesize DEL}}L(g, h)=  d^0\left[L\circ l_{g}+L\circ
r_{h}\circ i\right](\epsilon (x))_{| A_xG}
\]
with $\beta(g)=\alpha(h)=x$.

In conclusion, we have characterized the solutions of the
Lagrangian system determined by $L: G\longrightarrow \R$ as the
sequences
 $(g_1, \ldots, g_N)$, with $(g_k, g_{k+1})\in G_2$, for each $k\in\{1, \ldots, N-1\}$, and
\[
D_{\hbox{\footnotesize DEL}}L(g_k, g_{k+1})=0, \quad 1\leq k\leq
N-1.
\]

\subsection{Discrete Poincar\'e-Cartan sections}

Given a Lagrangian function $L:G\longrightarrow \R$, we will study
the geometrical properties of the discrete Euler-Lagrange
equations.

Consider the Lie algebroid $\pi^\tau:P^{\tau} G\cong
V\beta\oplus_G V\alpha \longrightarrow G$, and define the
\emph{Poincar\'e-Cartan 1-sections } $\Theta_L^-, \Theta_L^+\in
\Gamma ((\pi^\tau)^*)$ as follows
\begin{equation}\label{5.16'}
\Theta_L^-(g)(X_g, Y_g)= -X_g(L),\;\;\;\;\;  \Theta_L^+(g)(X_g,
Y_g)= Y_g(L),
\end{equation}
 for each $g\in G$ and $(X_g, Y_g)\in V_g\beta\oplus
V_g\alpha$. From the definition, we have that
\[\Theta_L^-(g)(X^{(1,0)}(g))=-\rvec{X}({g})(L)
\quand \Theta_L^-(g)(X^{(0,1)}(g))=0,
\]
and similarly
\[
\Theta_L^+(g)(X^{(0,1)}(g))=\lvec{X}({g})(L) \quand
\Theta_L^+(g)(X^{(1,0)}(g))=0,
\]
for $X\in\Gamma(\tau).$

We also have that $dL=\Theta_L^+ - \Theta_L^- $ and so, using
$d^2=0,$ it follows  that $d\Theta_L^+=d\Theta_L^-$. This means
that there exists a unique 2-section
$\Omega_L=-d\Theta_L^+=-d\Theta_L^-$, that will be called the
\emph{Poincar\'e-Cartan} 2-section. This 2-section will be
important for studying  symplecticity of the discrete
Euler-Lagrange equations.

\begin{proposition}\label{p5.2}
If $X$ and $Y$ are sections of the Lie algebroid $AG$ then
$$
\Omega_{L}(X^{(1,0)},Y^{(1,0)})=0\; , \ \Omega_{L}(X^{(0,1)},
Y^{(0,1)})=0, $$
 and
$$\Omega_{L}(X^{(1,0)},Y^{(0,1)})=-\rvec{X}(\lvec{Y}L) \ \hbox{ and   }\  \Omega_{L}(X^{(0,1)},Y^{(1,0)})=
{\rvec{Y}}(\lvec{X}L).
$$
\end{proposition}
\begin{proof}
A direct computation proves the result. \end{proof}

\begin{remark}\label{r4.2'}
{\bf Remark 4.3.} Let $g$ be an element of $G$ such that
$\alpha(g)=x$ and $\beta(g)=y$. Suppose that $U$ and $V$ are open
subsets of $M$, with $x\in U$ and $y\in V$, and that $\{X_i\}$ and
$\{Y_j\}$ are local bases of $\Gamma(\tau)$ on $U$ and $V$,
respectively. Then, $\{X_i^{(1,0)}, Y_j^{(0,1)}\}$ is a local
basis of $\Gamma(\pi^\tau)$ on the open subset $\alpha^{-1}(U)\cap
\beta^{-1}(V)$. Moreover, if we denote by
$\{(X^i)^{(1,0)},(Y^j)^{(0,1)}\}$ the dual basis of
$\{X_i^{(1,0)},Y_j^{(0,1)}\}$, we have that on the open subset
$\alpha^{-1}(U)\cap \beta^{-1}(V)$
\[
\begin{array}{l}
\Theta_L^-=-\rvec{X_i}(L)(X^i)^{(1,0)},\;\;\;\;\;
\Theta_L^+=\lvec{Y_j}(L)(Y^j)^{(0,1)},\\
\Omega_L=-\rvec{X_i}(\lvec{Y_j}(L))(X^i)^{(1,0)}\wedge
(Y^j)^{(0,1)}.
\end{array}
\]
\end{remark}

Finally, we obtain some useful expressions of the
Poincar\'e-Cartan 1-sections using the Lie algebroid isomorphisms
introduced in Theorems \ref{t3.3} and \ref{t3.4}.

We recall that the maps \[ \begin{array}{lcr}
A(\Phi^\beta):A(G_\beta*G_\alpha)\cong V\beta\oplus_G
V\alpha&\longrightarrow& A({\mathcal P}^\beta G)\cong {\mathcal P}^\beta(AG)\\
A(\Phi^\alpha):A(G_\beta*G_\alpha)\cong V\beta\oplus_G
V\alpha&\longrightarrow& A({\mathcal P}^\alpha G)\cong {\mathcal
P}^\alpha(AG)
\end{array}
\]
given by (\ref{e7}) and (\ref{3.13}) are Lie algebroid
isomorphisms over the identity of $G$. Moreover, if
$(v_{\epsilon(\beta(g))},Z_g)\in {\mathcal P}_g^\beta(AG)$ then,
from (\ref{e7}), it follows that
\begin{equation} \label{e10}
A_g(\Phi^{\beta})^{-1}(v_{\epsilon(\beta(g))}, Z_{g}) = (
Z_g-(T_{\epsilon(\beta(g))}l_{g})(v_{\epsilon(\beta(g))}),
(T_{\epsilon(\beta(g))}l_{g})(v_{\epsilon(\beta(g))})).
\end{equation}

On the other hand, if $(v_{\epsilon(\alpha(h))},Z_h)\in {\mathcal
P}_h^\alpha(AG)$ then, using (\ref{3.13}), we deduce that
\begin{equation} \label{e11}
\kern-12ptA_h(\Phi^{\alpha})^{-1}(v_{\epsilon(\alpha(h))}, Z_{h})
= (T_{\epsilon(\alpha(h))}(r_{h}\circ i)(v_{\epsilon(\alpha(h))}),
Z_h-T_{\epsilon(\alpha(h))}(r_{h}\circ
i)(v_{\epsilon(\alpha(h))})).
\end{equation}

Now, we introduce the sections $\Theta_L^\alpha\in
\Gamma((\tau^\alpha)^*)$ and $\Theta_L^\beta\in
\Gamma((\tau^\beta)^*)$ given by
\begin{equation}\label{5.18'}
\Theta_L^\alpha=(A(\Phi^\alpha)^{-1}, Id)^*(\Theta_L^{-}),\;\;\;
\Theta_L^\beta=(A(\Phi^\beta)^{-1}, Id)^*(\Theta_L^+).
\end{equation}

Using (\ref{e10}) and (\ref{e11}), we obtain that
\begin{align}\label{e11a}
\Theta_L^{\alpha}(h)(v_{\epsilon(\alpha(h))}, Z_{h})
&=-v_{\epsilon(\alpha(h))}(L\circ r_h\circ i),
\\
\label{e11b} \Theta_L^{\beta}(g)(v_{\epsilon(\beta(g))}, Z_{g})
&=v_{\epsilon(\beta(g))}(L\circ l_g),
\end{align}
for $(v_{\epsilon(\alpha(h))}, Z_h)\in {\mathcal P}_h^\alpha(AG)$
and $(v_{\epsilon(\beta(g))}, Z_g)\in {\mathcal P}_g^\beta(AG)$.

\subsubsection{Poincar\'e-Cartan $1$-sections: variational motivation}

Now, we follow a variational procedure to construct the
$1$-sections $\Theta_L^+$ and $\Theta_L^-$.
 We begin by calculating the extremals of ${\mathcal S} L$ for variations that do not fix
the point $g\in G$. For it, we consider the manifold
\[
{\mathcal C}^N=\{(g_1, \ldots, g_N)\in G^N\; / (g_k, g_{k+1})\in
G_2 \hbox{ for each } k, \;\; 1\leq k\leq N-1\}.
\]
If $c:(-\varepsilon,\varepsilon)\to {\mathcal C}^N$ is a curve in
${\mathcal C}^N$ and $c(0)=(g_1,\dots ,g_N)$ then there exist
$N+1$ curves $h_k:(-\varepsilon,\varepsilon)\to
\alpha^{-1}(\beta(g_k))$, for $ 0\leq k \leq N,$ with
$h_k(0)=\epsilon(\beta(g_k))$ and $g_0=g_1^{-1}$, such that
\[
c(t)=(h^{-1}_0(t)g_1h_1(t), h_1^{-1}(t)g_2 h_2(t), \ldots,
h_{N-2}^{-1}(t)g_{N-1} h_{N-1}(t), h_{N-1}^{-1}(t)g_Nh_N(t))
\]
for $t\in (-\varepsilon,\varepsilon).$ Thus, the tangent space to
${\mathcal C}^N$ at $(g_1,\dots ,g_N)$ may be identified with the
vector space $A_{\beta(g_0)}G\times A_{\beta(g_1)}G\times \dots
\times A_{\beta(g_N)}G,$ that is,
\[
T_{(g_1, g_2, \ldots, g_N)}{\mathcal C}^N\equiv\left\{(v_0, v_1,
\ldots, v_{N})\; /\; v_k\in A_{x_k}G, x_k=\beta(g_k), 0\leq k\leq
N\right\}.
\]
Now, proceeding as in Section \ref{section4.1}, we introduce the
action sum
\[{\mathcal S}L: {\mathcal C}^N\longrightarrow \R,\;\;\;\; SL(g_1,\dots ,g_N)=\sum_{k=1}^NL(g_k).
\]
Then,
\begin{eqnarray}\label{aaqq}
\frac{d}{dt}\Big|_{t=0}{\mathcal S} L(c(t))&=&
\sum_{k=1}^{N-1}\left[\d (L\circ l_{g_k})(\epsilon(x_k))(v_k)+\d (L\circ r_{g_{k+1}}\circ i)(\epsilon(x_{k+1}))(v_k)\right]\nonumber\\
&&+\d (L\circ r_{g_{1}}\circ i)(\epsilon(x_{0}))(v_0)+ \d (L\circ
l_{g_N})(\epsilon(x_N))(v_N).
\end{eqnarray}
Therefore, if $X_0,\dots ,X_N$ are sections of $\tau:AG\to M$
satisfying, $X_k(x_k)=v_k,$ for all $k,$ we have that
\begin{eqnarray*}
&&\frac{d}{dt}\Big|_{t=0}{\mathcal S}
L(c(t))=\sum_{k=1}^{N-1}\left[\lvec{X}_k({g_k})(L)-\rvec{X}_k({g_{k+1}})(L)
\right]-
\rvec{X_0}({g_1})(L)+\lvec{X}_N({g_N})(L)\\
&&\quad{}=\sum_{k=1}^{N-1}(D_{\hbox{\footnotesize DEL}}L (g_k,
g_{k+1}))(v_k) +\Theta_L^-(g_1)(X_0^{(1,0)}(g_1))+
 \Theta_L^+(g_N)(X_N^{(0,1)}(g_N)).
\end{eqnarray*}
Note that it is in the last two terms (that arise from the
boundary variations) where  appear the Poincar\'e-Cartan
1-sections.

\subsection{Discrete Lagrangian evolution operator}

We say that a differentiable mapping $\xi: G\longrightarrow G$ is
a \emph{discrete flow} or a \emph{discrete Lagrangian evolution
operator for $L$} if it verifies the following properties:
\begin{enumerate}
\item[-] $\hbox{graph}(\xi)\subseteq G_2$, that is, $(g, \xi(g))\in G_2$, $\forall g\in G$
($\xi$  is a second order operator).
\item[-] $(g, \xi(g))$ is a solution of the discrete Euler-Lagrange
equations, for all $g\in G$, that is,
$(D_{\hbox{DEL}}L)(g,\xi(g))=0,$ for all $g\in G.$
\end{enumerate}
In such a case
\begin{equation}\label{flow}
\d(L\circ l_{g}+L\circ r_{\xi(g)}\circ
i)(\epsilon(\beta(g)))_{|A_{\beta(g)}G}=0, \quad\text{for all
$g\in G$}
\end{equation}
or, in other terms,
\begin{equation}\label{5.22'}
\lvec{X}(g)(L)-\rvec{X}(\xi(g))(L)=0
\end{equation}
for every section $X$ of $AG$ and every $g\in G.$

Now, we define the prolongation ${\mathcal P}^{\tau}\xi:
V\beta\oplus_G V\alpha\longrightarrow V\beta\oplus_G V\alpha$   of
the second order operator $\xi: G\longrightarrow G$ as follows:
\begin{equation}
\label{4.9'} {\mathcal P}^{\tau}\xi=A(\Phi^{\alpha})^{-1}\circ
({Id}, T\xi)\circ A(\Phi^{\beta})
\end{equation}
with $A(\Phi^{\alpha})$ and $A(\Phi^{\beta})$ the isomorphisms
defined in Theorems \ref{t3.3} and \ref{t3.4} and
$(Id,T\xi):{\mathcal P}^\beta(AG)\to {\mathcal P}^\alpha(AG)$ the
map given by
\[
({Id}, T\xi)(v_{\epsilon(\beta(g))}, X_g)=
(v_{\epsilon(\beta(g))}, (T_g\xi)(X_g)), \mbox{ for }
(v_{\epsilon(\beta(g))}, X_g)\in {\mathcal P}_g^\beta(AG).
\]

Since the pair $(({Id}, T\xi), \xi)$ is a Lie algebroid morphism
between the Lie algebroids ${\mathcal
P}^{\beta}(AG)\longrightarrow G$ and ${\mathcal
P}^{\alpha}(AG)\longrightarrow G$ then the pair $({\mathcal
P}^{\tau}\xi, \xi)$ is also a Lie algebroid morphism
\[
\xymatrix{%
&V\beta\oplus_G V\alpha\ar[d]\ar[rr]^{\mathcal{P}^\tau\xi}&&V\beta\oplus_G V\alpha\ar[d]\\
&G\ar[rd]_{\beta}\ar[rr]^{\xi}&&G\ar[ld]^{\alpha}\\
&&M& }
\]
%
%\vspace{-0.5cm}
%
% \unitlength=1.00mm \special{em:linewidth 0.4pt}
%\linethickness{0.4pt}
%\begin{picture}(80.00,70.00)(-10,0)
%\put(40.00,30.00){\makebox(0,0)[cc]{$G$}}
%\put(44.00,30.33){\vector(1,0){31.00}}
%\put(80.00,30.00){\makebox(0,0)[cc]{$G$}}
%\put(42.00,26.33){\vector(1,-1){13.00}}
%\put(77.00,26.33){\vector(-1,-1){13.00}}
%\put(60.00,10.00){\makebox(0,0)[cc]{$M$}}
%\put(58.00,30.33){\makebox(0,0)[cb]{$\xi$}}
%\put(43.33,19.00){\makebox(0,0)[rc]{$\beta$}}
%\put(72.67,19.00){\makebox(0,0)[lc]{$\alpha$}}
%\put(40.00,55.00){\makebox(0,0)[cc]{$V\beta\oplus_G V\alpha$}}
%\put(80.00,55.00){\makebox(0,0)[cc]{$V\beta\oplus_G V\alpha$}}
%\put(52.00,55.00){\vector(1,0){15.33}}
%\put(40.00,50.00){\vector(0,-1){16.00}}
%\put(80.00,50.00){\vector(0,-1){16.00}}
%\put(58.00,56.00){\makebox(0,0)[cb]{${\mathcal P}^{\tau}\xi$}}
%\end{picture}
%\vspace{-.5cm}
From the definition of ${\mathcal P}^{\tau}\xi$, we deduce that
\begin{equation}\label{poi}
{\mathcal P}_g^{\tau}\xi(X_g, Y_g)=((T_g(r_{g\xi(g)}\circ
i))(Y_g), (T_g\xi)(X_g)+(T_g\xi)(Y_g)-T_g(r_{g\xi(g)}\circ
i)(Y_g))
\end{equation}
for all $(X_g, Y_g)\in V_g\beta\oplus V_g\alpha$. Moreover, from
(\ref{2.6'}) and  (\ref{poi}),
 we obtain that
\begin{equation}\label{zaa}
{\mathcal P}^{\tau}\xi(\rvec{X}(g),
\lvec{Y}(g))=(-\rvec{Y}(\xi(g)), (T_g\xi)(\rvec{X}(g)+\lvec{Y}(g))
+\rvec{Y}(\xi(g)))
\end{equation}
for all $X, Y$ sections of $AG$.

\subsection{Preservation of Poincar\'e-Cartan sections}

The following result  explains the sense in which the discrete
Lagrange evolution operator preserves the Poincar\'e-Cartan
2-section.

\begin{theorem}
\label{symplectic} Let $L: G\longrightarrow \R$ be a discrete
Lagrangian on a Lie groupoid $G$. Then:
\begin{enumerate}
\item The map $\xi$ is a discrete Lagrangian evolution operator for $L$ if and only if
$({\mathcal P}^{\tau}\xi,\xi)^*\Theta_L^-=\Theta_L^+$.
\item The map $\xi$ is a discrete Lagrangian evolution operator for $L$ if and only if $({\mathcal P}^{\tau}\xi,\xi)^*\Theta_L^--\Theta_L^-=dL$.
\item If $\xi$ is discrete Lagrangian evolution operator then $({\mathcal P}^{\tau}\xi,\xi)^*\Omega_L=\Omega_L$.
\end{enumerate}
\end{theorem}
\begin{proof}
From (\ref{5.18'}), it follows
\begin{equation}\label{5.24''*}(A(\Phi^{\alpha}),Id)^*(\Theta_L^{\alpha})=\Theta_L^-
,\;\;\;\;(A(\Phi^{\beta}),Id)^*(\Theta_L^{\beta})=\Theta_L^+.
\end{equation}
On the other hand, if $(v_{\epsilon(\beta(g))},X_g)\in {\mathcal
P}_g^\beta(AG)$ then, using (\ref{e11a}) and (\ref{e11b}), we have
that
$$\left\{((\hbox{Id}, T\xi), \xi)^*(\Theta_L^{\alpha})\right\}(g)(v_{\epsilon(\beta(g))}, X_g)=
-v_{\epsilon(\beta(g))}(L\circ r_{\xi(g)}\circ i)$$ and
$$\Theta^{\beta}_{L}(g)(v_{\epsilon(\beta(g))}, X_g)=v_{\epsilon(\beta(g))}(L\circ l_{g}).$$
Thus, $((\hbox{Id}, T\xi),
\xi)^*\Theta_L^{\alpha}=\Theta^{\beta}_{L}$  if and only if $\xi$
is a discrete Lagrangian evolution operator for $L.$ Therefore,
using this fact and (\ref{5.24''*}), we prove (i).

The second property follows from (i) by taking into account that
$dL=\Theta_L^+-\Theta_L^-$. Finally, (iii) follows using (ii) and
the fact that $({\mathcal P}^{\tau}\xi,\xi)$ is a Lie algebroid
morphism.

\end{proof}

\begin{remark}
Now, we present a  proof of the preservation of the
Poincar\'e-Cartan $2$-section using variational arguments. Given a
discrete Lagrangian evolution operator $\xi: G\longrightarrow G$
for $L,$ we may consider the function $\mathcal{S}_\xi L:G\to \R$
given by
\[
(\mathcal{S}_\xi L)(g)=L(g)+L(\xi(g)),\;\;\; \mbox{ for }g\in G.
\]
If $d$ is the differential on the Lie algebroid $V\beta\oplus_G
V\alpha\to G$ and $X,Y$ are sections of $\tau:AG\to M$ then, using
(\ref{5.16'}), (\ref{5.22'}) and (\ref{zaa}), we obtain that
\begin{eqnarray*}
d({{\mathcal S_\xi} L})(g)(\rvec{X}(g), \lvec{Y}(g))&=&
\rvec{X}(g)L+\lvec{Y}(g)L+(T_g\xi)(\rvec{X}(g))L+(T_g\xi)(\lvec{Y}(g))L\\
&=&
\rvec{X}(g)L+\lvec{Y}(g)L-\rvec{Y}(\xi(g))L+(T_g\xi)(\rvec{X}(g))L\\&&+(T_g \xi)(\lvec{Y}(g))L+\rvec{Y}(\xi(g))L\\
&=&-\Theta_L^-(g)(\rvec{X}(g), \lvec{Y}(g))+\left[({\mathcal
P}^{\tau}\xi,\xi)^*\Theta_L^+\right](g)(\rvec{X}(g), \lvec{Y}(g)).
\end{eqnarray*}
This implies that
\[({\mathcal P}^{\tau}\xi,\xi)^*\Theta_L^+-\Theta_L^- =d(\mathcal{S}_\xi L).
\]
Thus, we conclude that $({\mathcal
P}^{\tau}\xi,\xi)^*\Omega_L=\Omega_L$.
\end{remark}

\subsection{Lie groupoid morphisms and reduction}

Let  $(\Phi, \Phi_{0})$ be a  Lie groupoid morphism between the
Lie groupoids $G \rightrightarrows M$ and $G' \rightrightarrows
M'$. The prolongation $\mathcal{P}^\tau\Phi: V\beta\oplus_G
V\alpha\longrightarrow V\beta'\oplus_{G'} V\alpha'$ of the
morphism $(\Phi,\Phi_0)$ is defined by
\begin{equation}
\mathcal{P}^\tau_g\Phi (V,W)=(T_g\Phi(V), T_g\Phi(W))
\end{equation}
for every $(V,W)\in V_g\beta\oplus V_g\alpha$. It is easy to see
that $(\mathcal{P}^\tau\Phi,\Phi)$ is a morphism of Lie
algebroids.

\begin{theorem}
Let $(\Phi,\Phi_0)$ be a morphism of Lie groupoids from $G
\rightrightarrows M$ to $G'\rightrightarrows M'$. Let $L$ and $L'$
be discrete Lagrangian functions on $G$ and $G'$, respectively,
related by $L=L'\circ\Phi$. Then:
\begin{enumerate}
\item for every $(g,h)\in G_2$ and every $v\in A_{\beta(g)}G$ we have that
\begin{equation}
\label{DEL_tranformation}
D_{\mathrm{DEL}}L(g,h)(v)=D_{\mathrm{DEL}}L'\bigl(\Phi(g),\Phi(h)\bigr)\bigl(A_{\beta(g)}\Phi(v)\bigr).
\end{equation}
\item $(\mathcal{P}^\tau\Phi, \Phi)^*\Theta^{+}_{L'}=\Theta^{+}_{L},$
\item $(\mathcal{P}^\tau\Phi, \Phi)^*\Theta^{-}_{L'}=\Theta^{-}_{L},$
\item $(\mathcal{P}^\tau\Phi, \Phi)^*\Omega_{L'}=\Omega_{L}.$
\end{enumerate}
\end{theorem}
\begin{proof}
To prove the first we notice that, if $(\Phi,\Phi_0)$ is a
morphism of Lie groupoids, then we have that $\Phi\circ l_{g} =
l_{\Phi(g)}\circ\Phi$ and $\Phi\circ r_{h} =
r_{\Phi(h)}\circ\Phi$, from where we get
\begin{align*}
D_{\mathrm{DEL}}L(g,h)(v)
&=Tl_g(v)L+Tr_h(Ti(v))L\\
&=Tl_g(v)(L'\circ\Phi)+Tr_h(Ti(v))(L'\circ\Phi)\\
&=T\Phi(Tl_g(v))L'+T\Phi(Tr_h(Ti(v)))L'\\
&=Tl_{\Phi(g)}(T\Phi(v))L'+Tr_{\Phi(h)}(T\Phi(Ti(v)))L'\\
&=Tl_{\Phi(g)}(T\Phi(v))L'+Tr_{\Phi(h)}(Ti'(T\Phi(v)))L'\\
&=D_{\mathrm{DEL}}L'\bigl(\Phi(g),\Phi(h)\bigr)\bigl(A_{\beta(g)}\Phi(v)\bigr),
\end{align*}
where we have also used that $i'\circ\Phi=\Phi\circ i$ and
$A_{\beta(g)}\Phi(v)=T\Phi(v)$.

For the proof of the second, we have that
\begin{align*}
\pai{(\mathcal{P}^\tau\Phi, \Phi)^*\Theta^{+}_{L'}(g)}{(V,W)}
&=\pai{\Theta^{+}_{L'}(\Phi(g))}{(T_g\Phi(V),T_g\Phi(W))}\\
&=((T_g\Phi)(W))L' =WL =\pai{\Theta^{+}_{L}(g)}{(V,W)},
\end{align*}
for every $(V,W)\in\mathcal{P}^\tau_g G$. The proof of the third
is similar to the second, and finally, for the proof of (iv) we
just take the differential in (ii).
\end{proof}

As an immediate consequence of the above theorem we have that
\begin{corollary}\label{reduction}
Let $(\Phi,\Phi_0)$ be a morphism of Lie groupoids from
$G\rightrightarrows M$ to $G'\rightrightarrows M'$ and suppose
that $(g,h)\in G_2$.
\begin{enumerate}
\item If $(\Phi(g),\Phi(h))$ is a solution of the discrete Euler-Lagrange equations for $L'=L\circ\Phi$, then $(g,h)$ is a solution of the discrete Euler-Lagrange equations for $L$.
\item If $\Phi$ is a submersion then $(g,h)$ is a solution of the discrete Euler-Lagrange equations for $L$ if and only if $(\Phi(g), \Phi(h))$ is a solution of the discrete Euler-Lagrange equations for $L'$.
\item If $\Phi$ is an immersion, then $(g,h)$ is a solution of the discrete Euler-Lagrange equations for $L$ if and only if $D_{\mathrm{DEL}}L(\Phi(g),\Phi(h))$ vanishes over $\operatorname{Im}(A_\beta(g)\Phi)$.
\end{enumerate}
\end{corollary}
The case when $\Phi$ is an inmersion may be useful to modelize
holonomic mechanics on Lie groupoids, which is an imprescindible
tool for explicitely construct geometric integrators
(see~\cite{Hair,mawest}).

The particular case when $\Phi$ is a submersion is relevant for
reduction (see Section \ref{rdLs} in this paper).

\subsection{Discrete Legendre transformations}
\label{section5.6} Given a Lagrangian $L: G\longrightarrow \R$ we
define, just as the standard case \cite{mawest}, two
\emph{discrete Legendre transformations} $\F^{-}L:
G\longrightarrow A^*G$ and $\F^{+}L: G\longrightarrow A^*G$
 as follows
\begin{eqnarray*}
(\F^{-}L)(h)(v_{\epsilon(\alpha(h))})&=&-v_{\epsilon(\alpha(h))}(L\circ
r_h\circ i),\;\;\; \mbox{ for }
v_{\epsilon(\alpha(h))}\in A_{\alpha(h)}G,\\
(\F^{+}L)(g)(v_{\epsilon(\beta(g))})&=&v_{\epsilon(\beta(g))}(L\circ
l_g), \mbox{ for } v_{\epsilon(\beta(g))}\in A_{\beta(g)}G.
\end{eqnarray*}
\begin{remark}\label{r4.4'}
Note that $(\F^-L)(h)\in A^*_{\alpha(h)}G.$ Furthermore, if $U$ is
an open subset of $M$ such that $\alpha(h)\in U$ and $\{X_i\}$ is
a local basis of $\Gamma(\tau)$ on $U$ then
\[
\F^-L=\rvec{X_i}(L)(X^i\circ \alpha),\]
 on $\alpha^{-1}(U),$ where $\{X^i\}$ is the dual basis of
 $\{X_i\}$. In a similar way, if $V$ is an open subset of $M$ such
 that $\beta(g)\in V$ and $\{Y_j\}$ is a local basis of
 $\Gamma(\tau)$ on $V$ then
 \[
 \F^+L=\lvec{Y_j}(L)(Y^j\circ \beta),
 \]
 on $\beta^{-1}(V).$
\end{remark}
Next, we consider the prolongation $\tau^{\tau^*}: {\mathcal
P}^{\tau^*}(AG) \to A^*G$ of the Lie algebroid $\tau: AG \to M$
over the fibration $\tau^*: A^*G \to M$, that is,
\[
\begin{array}{rcl}
{\mathcal P}^{\tau^*}_{v^*}(AG)& =& \{ (v_{\tau^*(v^*)}, X_{v^*})
\in A_{\tau^*(v^*)} G \times T_{v^*}(A^*G) /
(T_{\tau^*(v^*)}\beta)(v_{\tau^*(v^*)}) \\&&=
(T_{v^*}\tau^*)(X_{v^*}) \}
\end{array}
\]
for $v^* \in A^*G$. Then, we may introduce the canonical section
$\Theta$ of the vector bundle $(\tau^{\tau^*})^*:({\mathcal
P}^{\tau^*}AG)^*\to A^*G$ as follows:
\begin{equation}\label{4.16'}
\Theta(v^*)(v_{\tau^*(v^*)}, X_{v^*})=v^*(v_{\tau^*(v^*)}),
\end{equation}
for $v^*\in A^*G$ and $(v_{\tau^*(v^*)}, X_{v^*})\in
\mathcal{P}^{\tau^*}_{v^*}(AG)$. $\Theta$ is called the
\emph{Liouville section}. Moreover, we define the \emph{canonical
symplectic section } $\Omega$ associated with $AG$ by
$\Omega=-d\Theta$, where $d$ is the differential on the Lie
algebroid $\tau^{\tau^*}:\mathcal{P}^{\tau^*}(AG)\to A^*G$. It is
easy to prove that $\Omega$ is nondegenerate and closed, that is,
it is a symplectic section of $\mathcal{P}^{\tau^*}(AG)$ (see
\cite{LMM,Medina}).

Now, let ${\mathcal P}^{\tau}\F^-L$ be the prolongation of $\F^-L$
defined by
\begin{equation}
\label{4.14'} {\mathcal P}^{\tau}\F^-L=(\hbox{Id}, T\F^-L)\circ
A(\Phi^{\alpha}): {\mathcal P}^{\tau}G\equiv V\beta\oplus_G
V\alpha\longrightarrow {\mathcal P}^{\tau^*}(AG),
\end{equation}
where $A(\Phi^\alpha):{\mathcal P}^\tau G\equiv V\beta \oplus_G
V\alpha\to {\mathcal  P}^\alpha(AG)$ is the Lie algebroid
isomorphism (over the identity of $G$) defined by (\ref{3.13}) and
$(Id,T\F^-L):{\mathcal P}^\alpha(AG)\to {\mathcal P}^{\tau^*}(AG)$
is the map given by
\[
(Id,T\F^{-}L)(v_{\epsilon(\alpha(h))},X_h)=(v_{\epsilon(\alpha(h))},
(T_h\F^-L)(X_h)),
\]
for $(v_{\epsilon(\alpha(h))},X_h)\in {\mathcal P}^\alpha_h(AG).$
Since the pair $((Id,T\F^-L),\F^-L)$ is a morphism between the Lie
algebroids ${\mathcal P}^\alpha(AG)\to G$ and ${\mathcal
P}^{\tau^*}(AG)\to A^*G$, we deduce that $({\mathcal
P}^\tau\F^-L,\F^-L)$ is also a morphism between the Lie algebroids
${\mathcal P}^\tau G\equiv V\beta\oplus_GV\alpha\to G$ and
${\mathcal P}^{\tau^*}(AG)\to A^*G.$ The following diagram
illustrates the above situation:
\[
\xymatrix{%
&V\beta\oplus_G V\alpha\ar[d]\ar[rr]^{\mathcal{P}^\tau\F^-L}&&\mathcal{P}^{\tau^*}(AG)\ar[d]\\
&G\ar[rd]_{\alpha}\ar[rr]^{\F^-L}&&A^*G\ar[ld]^{\tau^*}\\
&&M& }
\]
%
%\vspace{-.5cm}
%\begin{center}
%\unitlength=1.00mm \special{em:linewidth 0.4pt}
%\linethickness{0.4pt}
%\begin{picture}(80.00,70.00)(20,0)
%\put(40.00,30.00){\makebox(0,0)[cc]{$G$}}
%\put(44.00,30.33){\vector(1,0){31.00}}
%\put(80.00,30.00){\makebox(0,0)[cc]{$A^*G$}}
%\put(42.00,26.33){\vector(1,-1){13.00}}
%\put(77.00,26.33){\vector(-1,-1){13.00}}
%\put(60.00,10.00){\makebox(0,0)[cc]{$M$}}
%\put(58.00,31.33){\makebox(0,0)[cb]{$\F^{-}L$}}
%\put(43.33,19.00){\makebox(0,0)[rc]{$\alpha$}}
%\put(72.67,19.00){\makebox(0,0)[lc]{$\tau^*$}}
%\put(40.00,55.00){\makebox(0,0)[cc]{$V\beta\oplus_G V\alpha$}}
%\put(80.00,55.00){\makebox(0,0)[cc]{${\mathcal P}^{\tau^*}(AG)$}}
%\put(51.00,55.00){\vector(1,0){19.33}}
%\put(40.00,50.00){\vector(0,-1){16.00}}
%\put(80.00,50.00){\vector(0,-1){16.00}}
%\put(58.00,57.00){\makebox(0,0)[cb]{${\mathcal P}^{\tau}\F^-L$}}
%\end{picture}
%\end{center}
%\vspace{-1cm}

The prolongation ${\mathcal P}^{\tau}\F^-L$ can be explicitly
written as
\begin{equation}\label{eq19}
{\mathcal P}_h^{\tau}\F^-L(X_h, Y_h)=(T_h(i\circ r_{h^{-1}})(X_h),
(T_h\F^{-}L)(X_h)+(T_h\F^{-}L)(Y_h)),
\end{equation}
for $h\in G$ and $(X_h, Y_h)\in V_h\beta\oplus V_h\alpha$.

\begin{proposition}\label{p5.8'} If $\Theta$ is the Liouville
section of the vector bundle $({\mathcal P}^{\tau^*}(AG))^*\to
A^*G$ and $\Omega=-d\Theta$ is the canonical symplectic section of
$\wedge^2({\mathcal P}^{\tau^*}(AG))^*\to A^*G$ then
\begin{eqnarray*}
({\mathcal P}^{\tau}(\F^-L), \F^-L)^*\Theta= \Theta^-_L,&&
({\mathcal P}^{\tau}(\F^-L), \F^-L)^*\Omega= \Omega_L.
\end{eqnarray*}
\end{proposition}
\begin{proof}
Let $\Theta_L^\alpha$ be the section of $(\tau^\alpha)^*:
{\mathcal P}^\alpha (AG)^*\to G$ defined by (\ref{5.18'}). Then,
from (\ref{e11a}) and (\ref{4.16'}), we deduce that
\[
((Id,T\F^-L), \F^-L)^*\Theta=\Theta_L^\alpha.
\]
Thus, using (\ref{5.18'}), we obtain that
\[
({\mathcal P}^{\tau}\F^-L,\F^-L)^*\Theta=\Theta_L^-.
\]
Therefore, since the pair $({\mathcal P}^\tau\F^-L,\F^-L)$ is a
Lie algebroid morphism, it follows that
\[
(\mathcal{P}^\tau \F^{-}L,\F^{-}L)^*\Omega=\Omega_L.
\]
\end{proof}

 Now, we consider
the prolongation ${\mathcal P}^{\tau}\F^+L$ of $\F^+L$ defined by
\begin{equation}
\label{4.15'} {\mathcal P}^{\tau}\F^+L=(Id, T\F^+L)\circ
A(\Phi^{\beta}): {\mathcal P}^{\tau}G\equiv V\beta\oplus_G
V\alpha\longrightarrow {\mathcal P}^{\tau^*}(AG),
\end{equation}
where $A(\Phi^\beta):{\mathcal P}^\tau G\equiv V\beta\oplus_G
V\alpha\to {\mathcal P}^\beta(AG)$ is the Lie algebroid
isomorphism (over the identity of $G$) defined by (\ref{e7}) and
$(Id,T\F^+L):{\mathcal P}^\beta(AG)\to {\mathcal P}^{\tau^*}(AG)$
is the map given by
\[
(Id,T\F^+L)(v_{\epsilon(\beta(g))},X_g)=(v_{\epsilon(\beta(g))},(T_g\F^+L)(X_g)),
\]
for $(v_{\epsilon(\beta(g))},X_g)\in {\mathcal P}_g^\beta(AG).$ As
above, the pair $({\mathcal P}^\tau\F^+L,\F^+L)$ is a morphism
between the Lie algebroids ${\mathcal P}^\tau G\equiv
V\beta\oplus_G V\alpha\to G$ and ${\mathcal P}^{\tau^*}(AG)\to
A^*G$ and the following diagram illustrates the situation
\[
\xymatrix{%
&V\beta\oplus_G V\alpha\ar[d]\ar[rr]^{\mathcal{P}^\tau\F^+L}&&\mathcal{P}^{\tau^*}(AG)\ar[d]\\
&G\ar[rd]_{\beta}\ar[rr]^{\F^+L}&&A^*G\ar[ld]^{\tau^*}\\
&&M& }
\]
%\vspace{-1cm}
%\begin{center}
% \unitlength=1.00mm \special{em:linewidth 0.4pt}
%\linethickness{0.4pt}
%\begin{picture}(80.00,70.00)(30,0)
%\put(40.00,30.00){\makebox(0,0)[cc]{$G$}}
%\put(44.00,30.33){\vector(1,0){31.00}}
%\put(80.00,30.00){\makebox(0,0)[cc]{$A^*G$}}
%\put(42.00,26.33){\vector(1,-1){13.00}}
%\put(77.00,26.33){\vector(-1,-1){13.00}}
%\put(60.00,10.00){\makebox(0,0)[cc]{$M$}}
%\put(58.00,31.33){\makebox(0,0)[cb]{$\F^{+}L$}}
%\put(43.33,19.00){\makebox(0,0)[rc]{$\beta$}}
%\put(72.67,19.00){\makebox(0,0)[lc]{$\tau^*$}}
%\put(40.00,55.00){\makebox(0,0)[cc]{$V\beta\oplus_G V\alpha$}}
%\put(80.00,55.00){\makebox(0,0)[cc]{${\mathcal P}^{\tau^*}(AG)$}}
%\put(51.00,55.00){\vector(1,0){19.33}}
%\put(40.00,50.00){\vector(0,-1){16.00}}
%\put(80.00,50.00){\vector(0,-1){16.00}}
%\put(58.00,56.00){\makebox(0,0)[cb]{${\mathcal P}^{\tau}\F^+L$}}
%\end{picture}
%\end{center}
%\vspace{-1cm}

We also have:

\begin{proposition}\label{p5.8''} If $\Theta$ is the Liouville
section of the vector bundle $({\mathcal P}^{\tau^*}(AG))^*\to
A^*G$ and $\Omega=-d\Theta$ is the canonical symplectic section of
$\wedge^2({\mathcal P}^{\tau^*}(AG))^*\to A^*G$ then
\begin{eqnarray*}
({\mathcal P}^{\tau}(\F^+L), \F^+L)^*\Theta= \Theta^+_L,&&
({\mathcal P}^{\tau}(\F^+L), \F^+L)^*\Omega= \Omega_L.
\end{eqnarray*}
\end{proposition}
\begin{remark}
\label{evolution} $(i)$ If $\xi:G\to G$ is a smooth map then $\xi$
is a discrete Lagrangian evolution operator for $L$ if and only if
$\F^-L\circ\xi= \F^+L.$
\medskip

$(ii)$ If $(g,h)\in G_2$ we have that
\begin{equation}\label{5.32'}
(D_{\mbox{\footnotesize DEL}}L)(g,h)=\F^+L(g)-\F^-L(h).
\end{equation}
\end{remark}
%%%%%%%%%%%%%%%%%%%%%%%%%%%%%%%%%%%%%%%%%%%%%%

\subsection{Discrete regular Lagrangians}\label{regular}
First of all, we will introduce the notion of a discrete regular
Lagrangian.
\begin{definition}
 A Lagrangian $L:G\to \R$ on a Lie groupoid $G$ is said to be
regular if the Poincar{\'e}-Cartan $2$-section $\Omega_L$ is
symplectic on the Lie algebroid ${\mathcal P}^\tau G\equiv
V\beta\oplus_GV\alpha\to G.$
\end{definition}

Next, we will obtain necessary and sufficient conditions for a
discrete Lagrangian on a Lie groupoid to be regular.

\begin{theorem}\label{5.13}
Let $L:G\to \R$ be a Lagrangian function. Then:

{\bf a)} The following conditions are equivalent:
\begin{enumerate}
\item $L$ is regular.

\item The Legendre transformation $\F^-L$ is a local
diffeomorphism.

\item The Legendre transformation $\F^+L$ is a local
diffeomorphism.
\end{enumerate}

 {\bf b)} If $L:G\to \R$ is regular and $(g_0,h_0)\in
G_{2}$ is a solution of the discrete Euler-Lagrange equations for
$L$ then there exist two open subsets $U_0$ and $V_0$ of $G$, with
$g_0\in U_0$ and $h_0\in V_0,$ and there exists a (local) discrete
Lagrangian evolution operator $\xi_L:U_0\to V_0$ such that:
\begin{enumerate}
\item $\xi_L(g_0)=h_0,$

\item $\xi_L$ is a diffeomorphism and

\item $\xi_L$ is unique, that is, if $U_0'$ is an open subset of
$G$, with $g_0\in U_0'$ and $\xi_L':U'_0\to G$ is a (local)
discrete Lagrangian evolution operator then $ \xi'_{L|U_0\cap
U_0'}=\xi_{L|U_0\cap U_0'}$.
\end{enumerate}
\end{theorem}
\begin{proof}

\fbox{$a)$} First we will deduce the equivalence of the three
conditions

\fbox{(i)$\Rightarrow$(ii)} If $h\in G$, we need to prove that
$T_h(\F^-L): T_h G\longrightarrow T_{\F^-L(h)}A^*G$ is a linear
isomorphism. Assume that there exists $Y_h \in T_h G$ such that
$T_h(\F^-L)(Y_h)=0$. Since $\tau^*\circ \F^-L=\alpha,$ then
$(T_h\alpha)(Y_h)=0$, that is, $Y_h\in V_h\alpha$.

Therefore, $(0_h, Y_h)\in V_h\beta\oplus V_h\alpha$ and, from
(\ref{eq19}), we have that ${\mathcal P}^{\tau}_h(\F^- L)(0_h,
Y_h)=0$. Moreover, $({\mathcal
P}_h^{\tau}(\F^-L))^*\Omega(\F^-L(h))= \Omega_L(h)$ and
$\Omega(\F^-L(h))$ and $\Omega_L(h)$ are nondegenerate. Therefore,
we deduce that ${\mathcal P}^{\tau}_h(\F^- L)$ is a linear
isomorphism. This implies that $Y_h=0$. This proves that
$T_h(\F^-L):T_hG\to T_{\F^-L(h)}(A^*G)$ is a linear isomorphism.
In the same way we deduce \fbox{(i)$\Rightarrow$(iii)}.

\fbox{(ii)$\Rightarrow$(i)} We will assume that $\F^+L$ is a local
diffeomorphism, so that
\[
{\mathcal P}_g^{\tau}\F^+L: {\mathcal P}_g^{\tau}G\equiv
V_g\beta\oplus V_g\alpha\longrightarrow {\mathcal
P}_{\F^+L(g)}^{\tau^*}(AG)
\]
is a linear isomorphism, for all $g\in G$.

On the other hand, if $\Omega$ is the canonical symplectic section
of the vector bundle $\wedge^2({\mathcal P}^{\tau^*}(AG))^*\to
A^*G$ then, from Proposition \ref{p5.8''}, we deduce that
\[
({\mathcal P}_g^\tau\F^+L)^*(\Omega(\F^+L(g)))=\Omega_L(g).
\]
Thus, since $\Omega(\F^+L(g))$ is nondegenerate, we conclude that
$\Omega_L(g)$ is also nondegenerate, for all $g\in G.$ Using the
same arguments we deduce \fbox{(iii)$\Rightarrow$(i)}.

\fbox{$b)$}  Using Remark 4.11, we have that
\[
(\F^+L)(g_0)=(\F^-L)(h_0)=\mu_0\in A^*G.
\]
Thus, from the first part of this theorem, it follows that there
exit two open subsets $U_0$ and $V_0$ of $G$, with $g_0\in U_0$
and $h_0\in V_0$, and an open subset $W_0$ of $A^*G$ such that
$\mu_0\in W_0$ and
\[
\F^+L:U_0\to W_0,\;\;\; \F^-L:V_0\to W_0\] are diffeomorphisms.
Therefore, using Remark 4.11, we deduce that
\[
\xi_L=[(\F^-L)^{-1}\circ (\F^+L)]_{|U_0}:U_0\to V_0
\]
is a (local) discrete Lagrangian evolution operator. Moreover, it
is clear that $\xi_L(g_0)=h_0$ and, from the first part of this
theorem, we have that $\xi_L$ is a diffeomorphism.

Finally, if $U_0'$ is an open subset of $G$, with $g_0\in U_0'$,
and $\xi'_L: U_0'\to G$ is another (local) discrete Lagrangian
evolution operator then $\xi_{L|U_0\cap U_0'}':U_0\cap U_0'\to G$
is also a (local) discrete Lagrangian evolution operator.
Consequently, using Remark 4.11, we conclude that
\[
\xi_{L|U_0\cap U_0'}'=[(\F^-L)^{-1}\circ (\F^+L)]_{|U_0\cap
U_0'}=\xi_{L|U_0\cap U_0'}.
\]
\end{proof}

\begin{remark}\label{r4.9'}
{\rm  Using Remark 4.3, we deduce that the Lagrangian $L$ is
regular if and only if for every $g\in G$ and every local basis
$\{X_i\}$ (respectively, $\{Y_j\}$) of $\Gamma(\tau)$ on an open
subset $U$ (respectively, $V$) of $M$ such that $\alpha(g)\in U$
(respectively, $\beta(g)\in V$) we have that the matrix
$\rvec{X_i}(\lvec{Y_j}(L))$ is regular on $\alpha^{-1}(U)\cap
\beta^{-1}(V)$. }
\end{remark}

Let $L:G\to \R$ be a regular discrete Lagrangian on $G$. If
$f:G\to \R$ is a real $C^\infty$-function on $G$ then, using
Theorem \ref{5.13}, it follows that there exists a unique
$\xi_f\in \Gamma(\pi^\tau)$ such that
\begin{equation}\label{5.35}
i_{\xi_f}\Omega_L=df,
\end{equation}
$d$ being the differential of the Lie algebroid
$\pi^\tau:{\mathcal P}^\tau G\equiv V\beta\oplus_G V\alpha\to G$.
$\xi_f$ is called \emph{the Hamiltonian section associated to $f$
with respect to $\Omega_L$.}

Now, one may introduce a bracket of real functions on $G$ as
follows:

\begin{equation}\label{PoisL}
\{\cdot,\cdot\}_L:C^\infty(G)\times C^\infty(G)\to
C^\infty(G),\;\;\; \{f,g\}_L=-\Omega_L(\xi_f,\xi_g).
\end{equation}
Note that, from (\ref{5.35}) and Propositions \ref{p5.8'} and
\ref{p5.8''}, we obtain that
\begin{equation}\label{HamilL}
({\mathcal P}^\tau\F^{\pm}L)\circ \xi_{\bar{f}\circ \F^{\pm}L}=
\xi_{\bar{f}}\circ \F^{\pm}L,
\end{equation}
for $\bar{f}\in C^\infty(A^*G),$ where ${\mathcal
P}^\tau\F^{\pm}L:{\mathcal P}^\tau G\equiv
V\beta\oplus_GV\alpha\to {\mathcal P}^{\tau^*}(AG)$ is the
prolongation of $\F^{\pm}L$ (see Section \ref{section5.6}) and
$\xi_{\bar{f}}$ is the Hamiltonian section associated to the real
function $\bar{f}$ on $A^*G$ with respect to the canonical
symplectic section $\Omega$ on $\wedge^2({\mathcal
P}^{\tau^*}(AG))^*\to A^*G,$ that is,
$i_{\xi_{\bar{f}}}\Omega=d\bar{f}.$

On the other hand, we consider the canonical linear Poisson
bracket $\{\cdot,\cdot\}:C^\infty(A^*G)\times C^\infty(A^*G)\to
C^\infty(A^*G)$ on $A^*G$ defined by (see \cite{LMM})
\begin{equation}\label{clPb}
\{\bar{f}, \bar{g}\} = -\Omega (\xi_{\bar{f}}, \xi_{\bar{g}}), \;
\;  \mbox{ for } \bar{f}, \bar{g} \in C^{\infty}(A^*G).
\end{equation}
We have that (see \cite{LMM})
\[
\lcf \xi_{\bar{f}}, \xi_{\bar{g}} \rcf^{\tau^{\tau^*}} =
\xi_{\{\bar{f}, \bar{g}\}}.
\]
Moreover, from (\ref{PoisL}), (\ref{HamilL}), (\ref{clPb}) and
Propositions \ref{p5.8'} and \ref{p5.8''}, we deduce that
\[
\{\bar{f}\circ \F^{\pm}L,\bar{g}\circ
\F^{\pm}L\}_L=\{\bar{f},\bar{g}\}\circ \F^{\pm}L.
\]
Using the above facts, we may prove the following result.

\begin{proposition}\label{Poisson}
Let $L:G\to\R$ be a regular discrete Lagrangian.
\begin{enumerate}
\item
The Hamiltonian sections with respect to $\Omega_L$ form a Lie
subalgebra of the Lie algebra
$(\Gamma(\pi^\tau),\lcf\cdot,\cdot\rcf^{{\mathcal P}^\tau G}).$

\item The Lie groupoid $G$ endowed with the bracket $\{\cdot,
\cdot \}_L$ is a Poisson manifold, that is, $\{\cdot, \cdot\}_L$
is skew-symmetric, it is a derivation in each argument with
respect to the usual product of functions and it satisfies the
Jacobi identity.

\item The Legendre transformations $\F^{\pm}L:G\to A^*G$ are local
Poisson isomorphisms.
\end{enumerate}
\end{proposition}

\subsection{Discrete Hamiltonian evolution operator}

Let $L: G\to \R$ be a regular Lagrangian and assume, without the
loss of generality, that the Legendre transformations $\F^+L$ and
$\F^-L$ are global diffeomorphisms. Then,
$\xi_L=(\F^-L)^{-1}\circ(\F^+L)$ is the discrete Euler-Lagrange
evolution operator and one may define the \emph{discrete
Hamiltonian evolution operator}, $\tilde{\xi}_L: A^*G\to A^*G$, by
\begin{equation}\label{dheo}
\tilde{\xi}_L=\F^+L\circ \xi_L\circ (\F^+L)^{-1}\; .
\end{equation}

From Remark \ref{evolution}, we have  the following alternative
definitions
\[
\tilde{\xi}_L=\F^-L\circ \xi_L\circ (\F^-L)^{-1}, \qquad
\tilde{\xi}_L=\F^+L\circ (\F^-L)^{-1}
\]
of the discrete Hamiltonian evolution operator. The following
commutative diagram illustrates the situation

\unitlength=0.9mm \linethickness{0.4pt}
\begin{picture}(120.00,46.00)(0,20)
\put(55.00,60.00){\makebox(0,0)[cc]{$G$}}
\put(95.00,60.00){\makebox(0,0)[cc]{$G$}}
\put(58.67,60.00){\vector(1,0){30.33}}
\put(71.67,61.00){\makebox(0,0)[cb]{$\xi_{L}$}}
\put(31.00,30.00){\makebox(0,0)[cc]{$A^*G$}}
\put(75,30.00){\makebox(0,0)[cc]{$A^*G$}}
\put(119.00,30.00){\makebox(0,0)[cc]{$A^*G$}}
\put(52.00,56.67){\vector(-2,-3){16.00}}
\put(56.67,56.67){\vector(2,-3){16.00}}
\put(93.00,56.67){\vector(-2,-3){16.00}}
\put(98.33,56.67){\vector(2,-3){16.00}}
\put(37.00,29.67){\vector(1,0){31.33}}
\put(80.33,29.67){\vector(1,0){31.33}}
\put(40.33,45.00){\makebox(0,0)[rb]{$\F^-L$}}
\put(67.33,45.00){\makebox(0,0)[lb]{$\F^+L$}}
\put(89.00,45.00){\makebox(0,0)[lb]{$\F^-L$}}
\put(109.33,45.00){\makebox(0,0)[lb]{$\F^+L$}}
\put(50.00,33.13){\makebox(0,0)[cb]{$\tilde{\xi}_{L}$}}
\put(98.00,33.13){\makebox(0,0)[cb]{$\tilde{\xi}_{L}$}}
\end{picture}

Define the prolongation ${\mathcal
P}^{\tau^*}\tilde{\xi}_{L}:{\mathcal P}^{\tau^*}(AG)\to {\mathcal
P}^{\tau^*}(AG)$ of $\tilde{\xi}_{L}$  by
$${\mathcal
P}^{\tau^*}\tilde{\xi}_{L}={\mathcal P}^{\tau}\F^+{L}\circ
{\mathcal P}^{\tau}\xi_L\circ ({\mathcal P}^{\tau}\F^+{L})^{-1},$$
or, alternatively (see~\eqref{4.9'},~\eqref{4.14'}
and~\eqref{4.15'}),
\begin{equation}\label{pullopham}
{\mathcal P}^{\tau^*}\tilde{\xi}_{L}={\mathcal
P}^{\tau}\F^+{L}\circ ({\mathcal P}^{\tau}\F^-{L})^{-1}, \qquad
{\mathcal P}^{\tau^*}\tilde{\xi}_{L}={\mathcal
P}^{\tau}\F^-{L}\circ {\mathcal P}^{\tau}\xi_L\circ ({\mathcal
P}^{\tau}\F^-{L})^{-1}.
\end{equation}

%Observe that this map verifies
%\[
%{\mathcal P}^{\tau^*}\tilde{\xi}_{L} (v, X)=(v, T\tilde{\xi}_L(X))
%\]
%for all $(v, X)\in {\mathcal P}^{\tau^*}(AG)$.

\begin{proposition}
\label{p5.8.8} If $\Theta$ is the Liouville section of the vector
bundle $({\mathcal P}^{\tau^*}(AG))^*\to A^*G$ and
$\Omega=-d\Theta$ is the canonical symplectic section of
$\wedge^2({\mathcal P}^{\tau^*}(AG))^*\to A^*G$ then
\[
({\mathcal P}^{\tau^*}\tilde{\xi}_L, \tilde{\xi}_L)^*\Theta=
\Theta + d(L \circ (\F^-{L})^{-1}), \makebox[.4cm]{} ({\mathcal
P}^{\tau^*}\tilde{\xi}_L, \tilde{\xi}_L)^*\Omega= \Omega.
\]
%\begin{eqnarray*}
%({\mathcal P}^{\tau^*}\tilde{\xi}_L, \tilde{\xi}_L)^*\Theta=
%\Theta,&& ({\mathcal P}^{\tau^*}\tilde{\xi}_L,
%\tilde{\xi}_L)^*\Omega= \Omega.
%\end{eqnarray*}
Moreover, $\tilde{\xi}_L$ is a Poisson morphism for the canonical
Poisson bracket on $A^*G$.
\end{proposition}
\begin{proof}
The result follows using (\ref{dheo}), (\ref{pullopham}) and
Theorem \ref{symplectic} and Propositions \ref{p5.8'} and
\ref{Poisson}.
\end{proof}

\subsection{Noether's theorem}

Recall that classical Noether's theorem states that a continuous
symmetry of a Lagrangian leads to constants of the motion. In this
section, we prove  a discrete version of Noether's theorem, i.e.,
a theorem relating invariance of the discrete Lagrangian under
some transformation with the existence of  constants of the
motion.

\begin{definition}
A section $X$ of $AG$ is said to be a Noether's symmetry of the
Lagrangian $L$ if there exists a function $f\in\cinfty{M}$ such
that
\[
dL(X^{(1,1)})=\beta^*f-\alpha^*f.
\]
In this case, $L$ is said to be quasi-invariant under $X.$
\end{definition}

When  $dL(X^{(1,1)})=-\rvec{X}L + \lvec{X}L= 0$, we will say that
$L$ is \emph{invariant} under $X$ or that $X$ is an
\emph{infinitesimal symmetry} of the discrete Lagrangian $L$.

\begin{remark}
The infinitesimal invariance of the Lagrangian corresponds to a
finite invariance property as follows.
 Let $\Phi_s$ the flow of $\lvec{X}$ and
  $\gamma(s)=\Phi_s(\epsilon(x))$ be its integral curve  with $\gamma(0)=\epsilon(x),$ where $x=\beta(g)$.
  Then, the integral curve of $\lvec{X}$ at $g$ is
  $s\mapsto r_{\gamma(s)}g=g\gamma(s)$, and the integral curve of $-\rvec{X}$ through $\epsilon(x)$ is
  $s\mapsto \gamma(s)^{-1}$. On the other hand, if $(h,h')\in
  G_2$ and $Y_h\in V_h\beta,$ $Z_{h'}\in V_{h'}\alpha$ then
  \[
  (T_{(h,h')}m)(Y_h,Z_{h'})=(T_hr_{h'})(Y_h)+(T_{h'}l_h)(Z_{h'}).
  \]
  Using the above facts, we deduce that the integral curve $\mu$
  of the vector field $-\rvec{X} + \lvec{X}$ on $G$ satisfying
  $\mu(0)=g$ is
  \[
  \mu(s)=\gamma(s)^{-1}g\gamma(s), \mbox{ for all }s.
  \]
  Thus, the invariance of the Lagrangian may be written as
\[
L\bigl(\gamma(s)^{-1}g\gamma(s)\bigr)=L(g), \mbox{ for all } s.
\]
\end{remark}

If $L:G\to \R$ is a regular discrete Lagrangian, by a
\emph{constant of the motion } we mean a function $F$ invariant
under the discrete Euler-Lagrange evolution operator $\xi_L$, that
is, $F\circ\xi_L=F$.

\begin{theorem}[Discrete Noether's theorem]
If $X$ is a Noether symmetry of a discrete Lagrangian $L$, then
the function $F={\Theta^-_{L}}(X^{(1,1)})-\alpha^*f$ is a
\emph{constant of the motion } for the discrete dynamics defined
by $L$.
\end{theorem}
\begin{proof}
We first notice that ${\Theta^-_{L}}(X^{(1,1)})=\rvec{X}L$ so that
the function $F$ is $F=\rvec{X}L-\alpha^*f$.

If the Lagrangian $L$ is quasi-invariant under $X$ and $g$ is a
point in $G$, then
\[
-\rvec{X}(g)(L)+\lvec{X}(g)(L)=f(\beta(g))-f(\alpha(g)),
\]
so that
\[
\lvec{X}(g)(L)=\rvec{X}(g)(L)+f(\beta(g))-f(\alpha(g)).
\]
We substrate $\rvec{X}(\xi_L(g))(L)$ to both sides of the above
expression, so that
\begin{align*}
\lvec{X}(g)(L) - \rvec{X}(\xi_L(g))(L)
&=[\rvec{X}(g)(L)-f(\alpha(g))]-
[\rvec{X}(\xi(g))(L)-f(\alpha(\xi_L(g))]\\
&= F(g)-F(\xi_L(g)),
\end{align*}
from where the result immediately follows using (\ref{5.22'}).
\end{proof}

\begin{proposition}\label{5.8}
If $X$ is a Noether symmetry of the discrete Lagrangian $L$ then
\begin{equation}\label{5.24''}
{\mathcal L}_{X^{(1, 1)}}\Theta^-_{L}=d(\alpha^*f).
\end{equation}

Thus, if $L$ is regular, the complete lift $X^{(1, 1)}$ is a
Hamiltonian section with Hamiltonian function
$F=\Theta^-_{L}(X^{(1, 1)})-\alpha^*f$, i.e. $i_{X^{(1,
1)}}\Omega_{L}=dF$.
\end{proposition}
\begin{proof}
Indeed, if $dL (X^{(1, 1)})=\beta^*f-\alpha^*f$ and $Y$ is a
section of $AG$, we have that (see Proposition \ref{p5.2}),
\begin{eqnarray*}
({\mathcal L}_{X^{(1,1)}}\Theta_L^-)(Y^{(1,0)})&=&
-\Omega_L(X^{(1,1)},
Y^{(1,0)})+d( i_{X^{(1,1)}}\Theta_L^-)(Y^{(1,0)})\\
&=&-{\rvec{Y}}(\lvec{X}L)
+\rvec{Y}(\rvec{X}L)=\rvec{Y}(\alpha^*f-\beta^*f)\\
&=& d(\alpha^*f)(Y^{(1,0)}).
\end{eqnarray*}

On the other hand, using (\ref{RL}) and Proposition \ref{p5.2}, we
deduce that
\begin{eqnarray*} ({\mathcal
L}_{X^{(1,1)}}\Theta_L^-)(Y^{(0,1)})&=& -\Omega_L(X^{(1,1)},
Y^{(0,1)})+d(i_{X^{(1,1)}}\Theta_L^-)(Y^{(0,1)})\\
&=&-\rvec{X}(\lvec{Y}L)
+\lvec{Y}(\rvec{X}L)=[\lvec{Y},\rvec{X}](L)=0=
d(\alpha^*f)(Y^{(0,1)}).
\end{eqnarray*}
Thus, (\ref{5.24''}) holds. From (\ref{5.24''}), it follows that
\[
i_{X^{(1, 1)}}\Omega_{L}=-i_{X^{(1, 1)}}d\Theta^-_{L} =
di_{X^{(1,1)}}\Theta^-_{L}-{\mathcal L}_{X^{(1, 1)}}\Theta^-_{L} =
d[\Theta^-_{L}(X^{(1, 1)})-\alpha^*f]=dF,
\]
which completes the proof.
\end{proof}

We also have

\begin{proposition}\label{5.8'}
The vector space of Noether symmetries of the Lagrangian $L:G\to
\R$ is a Lie subalgebra of Lie algebra
$(\Gamma(\tau),\lcf\cdot,\cdot\rcf)$. \end{proposition}
\begin{proof}
Suppose that $X$ and $Y$ are Noether symmetries of $L$ and that
\begin{equation}\label{*1}
dL(X^{(1,1)})=-\rvec{X}L+\lvec{X}L=\beta^*f-\alpha^*f,
\end{equation}
\begin{equation}\label{*2}
dL(Y^{(1,1)})=-\rvec{Y}L+\lvec{Y}L=\beta^*g-\alpha^*g,
\end{equation}
with $f,g\in C^\infty (M)$. Then, using (\ref{RL}), (\ref{e1'})
and (\ref{e1''}), we have that
\begin{equation}\label{+}
dL(\lcf
X,Y\rcf^{(1,1)})=\rvec{X}(\rvec{Y}L)-\rvec{Y}(\rvec{X}L)+\lvec{X}(\lvec{Y}L)-\lvec{Y}(\lvec{X}L).
\end{equation}

On the other hand, from (\ref{linv}), (\ref{rinv}), (\ref{*1}) and
(\ref{*2}), we deduce that
\[\begin{array}{lll}
\rvec{X}(\rvec{Y}L)= \rvec{X}(\lvec{Y}L)-\alpha^*(\rho(X)(g)),&&
\rvec{Y}(\rvec{X}L)=\rvec{Y}(\lvec{X}L)-\alpha^*(\rho(Y)(f)), \\
\lvec{X}(\lvec{Y}L)= \lvec{X}(\rvec{Y}L)+\beta^*(\rho(X)(g)),&&
\lvec{Y}(\lvec{X}L)=\lvec{Y}(\rvec{X}L)+\beta^*(\rho(Y)(f)).
\end{array}
\]
Thus, using (\ref{RL}) and (\ref{+}), we obtain that
\[
dL(\lcf
X,Y\rcf^{(1,1)})=\beta^*(\rho(X)(g)-\rho(Y)(f))-\alpha^*(\rho(X)(g)-\rho(Y)(f)).
\]
Therefore, $\lcf X, Y\rcf$  is a Noether symmetry of $L.$
\end{proof}

\begin{remark}
If $L:G\to \R$ is a regular discrete Lagrangian then, from
Propositions \ref{5.8} and \ref{5.8'}, it follows that the
complete lifts of Noether symmetries of $L$ are a Lie subalgebra
of the Lie algebra of Hamiltonian sections with respect to
$\Omega_L$.
\end{remark}

\section{Examples}
\subsection{Pair or Banal groupoid}

We consider the pair (banal) groupoid $G=M\times M$, where the
structural maps are
\begin{eqnarray*}
&\alpha(x, y)=x, \; \; \beta(x, y)=y, \; \; \epsilon(x)=(x,x), \; \; i(x, y)=(y, x),&\\
&m((x, y), (y, z))=(x, z).&
\end{eqnarray*}

We know that the Lie algebroid of $G$ is isomorphic to the
standard Lie algebroid $\tau_{M}: TM \to M$ and the map
\[
\Psi: AG = V_{\epsilon(M)}\alpha \to TM, \makebox[.4cm]{} (0_{x},
v_{x}) \in T_{x}M \times T_{x}M \to \Psi_{x}(0_{x}, v_{x}) =
v_{x},\mbox{ for $x \in M$,}
\]
 induces an isomorphism (over the
identity of $M$) between $AG$ and $TM$. If $X$ is a section of
$\tau_{M}: AG \simeq TM \to M$, that is, $X$ is a vector field on
$M$ then $\rvec{X}$ and $\lvec{X}$ are the vector fields on $M
\times M$ given by
\[
\rvec{X}(x, y)=(-X(x), 0_y) \in T_{x}M \times T_{y}M \quand
\lvec{X}(x, y)=(0_x, X(y)) \in T_{x}M \times T_{y}M,
\]
for $(x, y) \in M \times M$. On the other hand, if $(x, y) \in M
\times M$ we have that the map
\[
\begin{matrix}
&{\mathcal P}^{\tau_{M}}_{(x, y)}G \equiv V_{(x, y)}\beta \oplus
V_{(x, y)}\alpha &\to &T_{(x, y)}(M \times M) \simeq T_{x}M \times
T_{y}M,\\
&((v_{x}, 0_{y}), (0_{x}, v_{y})) &\to &(v_{x}, v_{y})
\end{matrix}
\]
induces an isomorphism (over the identity of $M \times M$) between
the Lie algebroids $\pi^{\tau_{M}}: {\mathcal P}^{\tau_{M}}G
\equiv V\beta \oplus_{G} V\alpha \to G = M \times M$ and $\tau_{(M
\times M)}: T(M \times M) \to M \times M$.

Now, given a discrete Lagrangian $L: M \times M \to \R$ then the
discrete Euler-lagrange equations for $L$ are:
\begin{equation}\label{5.0*}
\lvec{X}(x, y)(L) - \rvec{X}(y, z)(L)=0,\;  \hbox{ for all } X\in
\mathfrak{X}(M),
\end{equation}
which are equivalent to the classical discrete Euler-Lagrange
equations
\[
D_2L(x, y)+D_1L(y, z)=0
\]
(see, for instance, \cite{mawest}). The Poincar\'e-Cartan
$1$-sections $\Theta_{L}^{-}$ and $\Theta_{L}^{+}$ on
$\pi^{\tau_{M}}: {\mathcal P}^{\tau_{M}}G \simeq T(M \times M) \to
G = M \times M$ are the $1$-forms on $M \times M$ defined by
\[
\Theta_{L}^{-}(x, y)(v_{x}, v_{y}) = -v_{x}(L), \makebox[.3cm]{}
\Theta_{L}^{+}(x, y)(v_{x}, v_{y}) = v_{y}(L),
\]
for $(x, y) \in M \times M$ and $(v_{x}, v_{y}) \in T_{x}M \times
T_{y}M \simeq T_{(x, y)}(M \times M)$.

In addition, if $\xi: G = M \times M \to G=M \times M$ is a
discrete Lagrangian evolution operator then the prolongation of
$\xi$
\[
{\mathcal P}^{\tau_{M}}\xi: {\mathcal P}^{\tau_{M}}G \simeq T(M
\times M) \to {\mathcal P}^{\tau_{M}}G \simeq T(M \times M)
\]
is just the tangent map to $\xi$ and, thus, we have that
\[
\xi^*\Omega_L=\Omega_L,
\]
$\Omega_{L} = -d\Theta_{L}^{-} = -d\Theta_{L}^{+}$ being the
Poincar\'{e}-Cartan $2$-form on $M \times M$. The Legendre
transformations $\F^{-}L: G = M \times M \to A^{*}G \simeq T^*M$
and $\F^{+}L: G = M \times M \to A^{*}G \simeq T^*M$ associated
with $L$ are the maps given by
\[
\F^-L(x, y) = -D_{1}L (x, y) \in T_{x}^{*}M, \makebox[.4cm]{}
\F^+L(x, y) = D_{2}L (x, y) \in T_{y}^{*}M
\]
for $(x, y) \in M \times M$. The Lagrangian $L$ is regular if and
only if the matrix $\displaystyle{\left( \frac{\partial^2
L}{\partial x\partial y}\right)}$ is regular. Finally, a Noether
symmetry is a vector field $X$ on $M$ such that
\[
D_{1}L(x, y)(X(x)) + D_{2}L(x, y)(X(y)) = f(y) - f(x),
\]
for $(x, y) \in M \times M$, where $f: M \to \R$ is a real
$C^{\infty}$-function on $M$. If $X$ is a Noether symmetry then
\[
x \to F(x)  = D_{1}L(x, y)(X(x)) - f(x)
\]
is a constant of the motion.

In conclusion, we recover all the geometrical formulation of the
classical discrete Mechanics on the discrete state space $M \times
M$ (see, for instance, \cite{mawest}).

\subsection{Lie groups}
We consider a Lie group $G$ as a groupoid  over one point
$M=\{\frak e \}$, the identity element of $G$. The structural maps
are
\[
\alpha(g)={\frak e} , \; \; \beta(g)={\frak e} , \; \;
\epsilon({\frak e})={\frak e}, \; \; i(g)=g^{-1}, \; \; m(g,
h)=gh, \; \; \mbox{ for } g, h \in G.
\]
The Lie algebroid associated with $G$ is just the Lie algebra
${\frak g}=T_{\frak e}G$ of $G$. Given $\xi\in {\frak g}$ we have
the left and right invariant vector fields:
\[
\lvec{\xi}(g)=(T_{\frak e}l_g)(\xi),\ \ \rvec{\xi}(g)=(T_{\frak e}
r_g)(\xi), \; \; \mbox{ for } g \in G.
\]
Thus, given a Lagrangian $L: G\longrightarrow \R$ its discrete
Euler-Lagrange equations are:
\[
(T_{\frak e}l_{g_k})(\xi)(L)-(T_{\frak e}r_{g_{k+1}})(\xi)(L)=0,\;
\hbox{ for all } \xi\in {\frak g} \mbox{ and } g_{k}, g_{k+1} \in
G,
\]
or, $(l_{g_k}^*dL)({\frak e}) = (r_{g_{k+1}}^*dL)({\frak e})$.
Denote by $\mu_k=(r_{g_k}^*dL)({\frak e})$ then the discrete
Euler-Lagrange equations are written as
\begin{equation}\label{5.0'*}
\mu_{k+1}=Ad^*_{g_{k}}\mu_k,
\end{equation}
where $\map{Ad}{G\times\mathfrak{g}}{\mathfrak{g}}$ is the adjoint
action of $G$ on $\mathfrak{g}$. These equations are known as the
\emph{discrete Lie-Poisson equations } (see
\cite{BoSu,MPS1,MPS2}).

Finally, an infinitesimal symmetry of $L$ is an element $\xi\in
\mathfrak{g}$ such that $(T_{\frak e}l_{g})(\xi)(L)=(T_{\frak
e}r_{g})(\xi)(L)$, and then the associated constant of the motion
is $F(g)=(T_{\frak e}l_{g})(\xi)(L)=(T_{\frak e}r_{g})(\xi)(L)$.
Observe that all the Noether's symmetries are infinitesimal
symmetries of $L$.

\subsection{Transformation or action Lie groupoid}
\label{TaLg}

Let $H$ be a Lie group and $\cdot:M\times H\to M$, $(x,h)\in
M\times H\mapsto xh,$ a right action of $H$  on $M$. As we know,
$H$ is a Lie groupoid over the identity element ${\frak e}$  of
$H$ and we will denote by $\alpha,\beta, \epsilon, m$ and $i$ the
structural maps of $H$. If $\pi:M\to \{{\frak e}\}$ is the
constant map then is clear that the space
\[
M_{\; \pi}\times_\alpha H=\{(x,h)\in M\times H/\pi(x)=\alpha(h)\}
\]
is the cartesian product $G=M\times H$ and that $\cdot :M\times
H\to M$ induces an action of the Lie groupoid $H$ over the map
$\pi:M\to\{{\frak e}\}$ in the sense of Section \ref{section2.2}
(see Example $6$ in Section \ref{section2.2}). Thus, we may
consider the action Lie groupoid $G=M\times H$ over $M$ with
structural maps given by
\begin{equation}\label{*tilde}
\begin{array}{l}
\tilde{\alpha}_\pi(x,h)=x,\;\;\; \tilde{\beta}_\pi(x,h)=xh,\;\;\;
\tilde{\epsilon }_\pi(x)=(x,{\frak e}),\\
\tilde{m}_\pi((x,h),(xh,h'))=(x,hh'),\;\;\;
\tilde{i}_\pi(x,h)=(xh, h^{-1}).
\end{array}
\end{equation}
Now, let ${\mathfrak h}=T_{{\frak e}}H$ be the Lie algebra of $H$
and $\Phi:{\mathfrak h}\to {\mathfrak X}(M)$ the map given by
\[
\Phi(\eta)=\eta_M,\;\;\;\; \mbox{for } \eta\in {\mathfrak h},
\]
where $\eta_M$ is the infinitesimal generator of the action
$\cdot:M\times H\to M$ corresponding to $\eta$. Then, $\Phi$
defines an action of the Lie algebroid ${\mathfrak{ h}}\to \{
\mbox{a point}\}$ over the projection $\pi:M\to \{\mbox{a
point}\}$ and the corresponding action Lie algebroid $pr_1:M\times
{\mathfrak{ h}}\to M$ is just the Lie algebroid of $G=M\times H$
(see Example $6$ in Section \ref{section2.2}).

We have that $\Gamma(pr_1)\cong \{\tilde{\eta}:M\to {\mathfrak
h}/\tilde{\eta} \mbox{ is  smooth }\}$ and that the Lie algebroid
structure $(\lcf\cdot,\cdot\rcf_\Phi,\rho_\Phi)$ on $pr_1:M\times
H\to M$ is given by
\[
\lcf
\tilde{\eta},\tilde{\mu}\rcf_{\Phi}(x)=[\tilde{\eta}(x),\tilde{\mu}(x)]
+
(\tilde{\eta}(x))_M(x)(\tilde{\mu})-(\tilde{\mu}(x))_M(x)(\tilde{\eta}),\;\;\;
\rho_\Phi(\tilde{\eta})(x)=(\tilde{\eta}(x))_M(x),
\]
for $\tilde{\eta}, \tilde{\mu}\in \Gamma(pr_1)$ and $x\in M.$
Here, $[\cdot,\cdot]$ denotes the Lie bracket of ${\frak h}$.

If $(x,h)\in G=M\times H$ then the left-translation
$l_{(x,h)}:\tilde{\alpha}_\pi^{-1}(xh)\to
\tilde{\alpha}_\pi^{-1}(x)$ and the right-translation
$r_{(x,h)}:\tilde{\beta}_\pi^{-1}(x)\to
\tilde{\beta}_\pi^{-1}(xh)$ are given
\begin{equation}\label{tilde+}
l_{(x,h)}(xh,h')=(x,hh'),\;\;\;
r_{(x,h)}(x(h')^{-1},h')=(x(h')^{-1},h'h).
\end{equation}

Now, if $\eta\in {\mathfrak h}$ then $\eta$ defines a constant
section $C_\eta:M\to {\mathfrak h}$ of $pr_1:M\times {\mathfrak
h}\to M$ and, using (\ref{linv}), (\ref{rinv}), (\ref{*tilde}) and
(\ref{tilde+}), we have that the left-invariant and the
right-invariant vector fields $\lvec{C}_\eta$ and $\rvec{C}_\eta$,
respectively, on $M\times H$ are defined by
\begin{equation}\label{*+}
\rvec{C}_\eta(x,h)=(-\eta_M(x),\rvec{\eta}(h)),\;\;\;\;
\lvec{C}_\eta(x,h)=(0_x,\lvec{\eta}(h)),
\end{equation}
for $(x,h)\in G=M\times H.$

Note that if $\{\eta_i\}$ is a basis of ${\mathfrak h}$ then
$\{C_{\eta_i}\}$ is a global basis of $\Gamma(pr_1).$

Next, suppose that $L:G=M\times H\to \R$ is a Lagrangian function
and for every $h\in H$ (resp., $x\in M$) we will denote by $L_h$
(resp., $L_x$) the real function on $M$ (resp., on $H$) given by
$L_h(y)=L(y,h)$ (resp., $L_x(h')=L(x,h'))$. Then, a composable
pair $((x,h_k),(xh_k,h_{k+1}))\in G_2$ is a solution of the
discrete Euler-Lagrange equations for $L$ if
\[
\lvec{C}_\eta(x,h_k)(L)-\rvec{C}_\eta(xh_k,h_{k+1})(L)=0, \mbox{
for all } \eta\in {\mathfrak h},
\]
or, in other terms (see (\ref{*+}))
\[
\{(T_{{\frak e}} l_{h_k})(\eta)\}(L_x)-\{(T_{{\frak e}}
r_{h_{k+1}})(\eta)\}(L_{xh_k})+\eta_M(xh_k)(L_{h_{k+1}})=0, \mbox
{ for all } \eta\in {\mathfrak h.}
\]

As in the case of Lie groups, denote by $\mu_k(x,h_k)=d(L_x\circ
r_{h_k})({\frak e}).$ Then, the discrete Euler-Lagrange equations
for $L$ are written as
\[
\mu_{k+1}(xh_k,h_{k+1})=Ad^*_{h_k}\mu_k(x,h_k)+d(L_{h_{k+1}}\circ
((x{h_k})\cdot))(e),
\]
where $(xh_k)\cdot:H\to M$ is the map defined by
\[
(xh_k)\cdot (h)=x(h_kh),\mbox{ for } h\in H.
\]
In the particular case when $M$ is the orbit of $a\in V$ under a
representation of $G$ on a real vector space $V$, the resultant
equations were obtained by Bobenko and Suris,
see~\cite{BoSu,BoSu2}, and they were called the \emph{discrete
Euler-Poincar\'e equations}.

Finally, an element  $\xi\in\mathfrak{h}$ is an infinitesimal
symmetry of $L$ if
\[
\xi_M(x)(L_{h})-\rvec{\xi}(h)(L_x) + \lvec{\xi}(h)(L_x)=f(xh)-f(x)
\]
where $f: M\longrightarrow \R$ is a real $C^\infty$-function on
$M$. The associated constant of the motion is
\[
F(x, h)=-\xi_M(x)(L_{h})+\rvec{\xi}(h)(L_x) -f(x),
\]
for $(x,h)\in M\times H$.

\subsubsection*{The heavy top}
As a concrete example of a system on a transformation Lie groupoid
we consider a discretization of the heavy top. In the continuous
theory~\cite{mart}, the configuration manifold is the
transformation Lie algebroid $\tau:S^2\times\mathfrak{so}(3)\to
S^2$ with Lagrangian
\[
L_c(\Gamma,\Omega)=\frac{1}{2}\Omega\cdot
I\Omega-mgl\Gamma\cdot\e,
\]
where $\Omega\in\R^3\simeq\mathfrak{so}(3)$ is the angular
velocity, $\Gamma$ is the direction opposite to the gravity and
$\e$ is a unit vector in the direction from the fixed point to the
center of mass, all them expressed in a frame fixed to the body.
The constants $m$, $g$ and $l$ are respectively the mass of the
body, the strength of the gravitational acceleration and the
distance from the fixed point to the center of mass. The matrix
$I$ is the inertia tensor of the body. In order to discretize this
Lagrangian it is better to express it in terms of the matrices
$\hat{\Omega}\in\mathfrak{so}(3)$ such that
$\hat{\Omega}v=\Omega\times v$. Then
\[
L_c(\Gamma,\Omega)=\frac{1}{2}\tr(\hat{\Omega}\I\hat{\Omega}^T)-mgl\Gamma\cdot\e.
\]
where $\I=\frac{1}{2}\tr(I) I_3-I$. We can define a discrete
Lagrangian $L: G = S^2\times SO(3)\to\R$ for the heavy top  by
\[
L(\Gamma_k,W_k)=-\frac{1}{h}\tr(\I W_k)-hmgl\Gamma_k\cdot\e.
\]
which is obtained by the rule
$\hat{\Omega}=R^T\dot{R}\approx\frac{1}{h}R_{k}^T(R_{k+1}-R_k)=\frac{1}{h}(W_k-I_3)$,
where $W_k=R_k^TR_{k+1}$.

The value of the action on an admissible variation is
\begin{align*}
\lambda(t)&=L(\Gamma_k,W_ke^{tK})+L(e^{-tK}\Gamma_{k+1},e^{-tK}W_{k+1})\\
&=-\frac{1}{h}\left[ \tr(\I
W_ke^{tK})+mglh^2\Gamma_k\cdot\e+\tr(\I e^{-tK}W_{k+1}) +mglh^2
(e^{-tK}\Gamma_{k+1})\cdot\e\right],
\end{align*}
where $\Gamma_{k+1}=W^T\Gamma_{k}$ (since the above pairs must be
composable) and $K\in\mathfrak{so}(3)$ is arbitrary. Taking the
derivative at $t=0$ and after some straightforward manipulations
we get the DEL equations
\[
M_{k+1}-W_k^T M_kW_k +mglh^2(\widehat{\Gamma_{k+1}\times\e})=0
\]
where $M=W\I-\I W^T$. Finally, in terms of the axial vector $\Pi$
in $\R^3$ defined by $\hat{\Pi}=M$, we can write the equations in
the form
\[
\Pi_{k+1}=W_k^T\Pi_k+mglh^2\Gamma_{k+1}\times\e.
\]
\begin{remark}
The above equations are to be solved as follows. From
$\Gamma_k,W_k$ we obtain $\Gamma_{k+1}=W_k\Gamma_k$ and $\Pi_k$
from $\hat{\Pi}_k=W_k\I-\I W_k^T$. The DEL equation gives
$\Pi_{k+1}$ in terms of the above data. Finally we get $W_{k+1}$
as the solution of the equation $\hat{\Pi}_{k+1}=W_{k+1}\I-\I
W_{k+1}^T$, as in~\cite{Mose}.
\end{remark}

In the continuous theory, the section $X(\Gamma)=(\Gamma,\Gamma)$
of $S^2\times\mathfrak{so}(3)\to S^2$ is a symmetry of the
Lagrangian (see~\cite{mart}). We will show next that such a
section is also a symmetry of the discrete Lagrangian. Indeed, it
is easy to see that the left and right vector fields associated to
$X$ coincide $\rvec{X}=\lvec{X}$ and are both equal to
\[
\rvec{X}(\Gamma,W)=\Bigl((\Gamma,0),(W,\hat{\Gamma}W))\Bigr)\in TG
=  TS^2\times TSO(3).
\]
Thus $\rho^{{\mathcal P}^{\tau}G}(X^{(1,1)})=0$ so that $X$ is a
symmetry of the Lagrangian. In fact it is a symmetry of any
discrete Lagrangian defined on $G = S^2\times SO(3)$. The
associated constant of motion is
\[
(\rvec{X}L)(W,\Gamma) =\tr(\I\hat{\Gamma}W)
=\frac{1}{2}\tr[(W\I-\I W^T)\hat{\Gamma}] =-\Pi\cdot\Gamma,
\]
i.e. (minus) the angular momentum in the direction of the vector
$\Gamma$.

\subsection{Atiyah or gauge groupoids}

Let $p: Q \to M$ be a principal $G$-bundle. A \emph{discrete
connection} on $p:Q\to M$ is a map ${\mathcal A}_{d}: Q \times Q
\to G$ such that
\begin{equation}
\label{Ad} {\mathcal A}_{d}(gq, hq') = h{\mathcal A}_{d}(q,
q')g^{-1} \qquand \mathcal{A}_d(q,q)=\mathfrak{e}
\end{equation}
for $g, h \in G$ and $q, q' \in Q$, $\mathfrak{e}$ being the
identity in the group $G$ (see \cite{TL,LeMaWe}). We remark that a
discrete principal connection may be considered as the discrete
version of an standard (continuous) connection on $p: Q \to M$. In
fact, if ${\mathcal A}_{d}:Q\times Q\to G$ is such a connection
then it induces, in a natural way, a continuous connection
${\mathcal A}_{c}: TQ \to {\mathfrak g}$ defined by
\[
{\mathcal A}_{c}(v_{q}) = (T_{(q, q)}{\mathcal A}_{d})(0_{q},
v_{q}),
\]
for $v_{q} \in T_{q}Q$. Moreover, if we choose a local
trivialization of the principal bundle $p: Q \to M$ to be $G
\times U$, where $U$ is an open subset of $M$ then, from
(\ref{Ad}), it follows that there exists a map $A: U \times U \to
G$ such that
\[
{\mathcal A}_{d}((g, x), (g', y)) = g'A(x, y)g^{-1}, \qquand A(x,
x) = {\frak e},
\]
for $(g, x), (g', x') \in G \times U$ (for more details, see
\cite{TL,LeMaWe}).

On the other hand, using the discrete connection ${\mathcal
A}_{d}$, one may identify the open subset $(p^{-1}(U) \times
p^{-1}(U))/G \simeq ((G \times U)\times (G \times U))/G$ of the
Atiyah groupoid $(Q \times Q)/G$ with the product manifold $(U
\times U) \times G$. Indeed, it is easy to prove that the map
\[
((G \times U)\times (G \times U))/G \to (U \times U)\times G,
\]
\[
[((g, x), (g', y))] \to ((x, y), {\mathcal A}_{d}((e, x),
(g^{-1}g', y))) = ((x, y), g^{-1}g'A(x, y)),
\]
is bijective. Thus, the restriction to $((G \times U)\times (G
\times U))/G$ of the Lie groupoid structure on $(Q \times Q)/G$
induces a Lie groupoid structure in $(U \times U) \times G$ with
source, target and identity section given by
\[
\begin{array}{lr}
\alpha: (U \times U) \times G \to U; & ((x, y), g) \to x,
\\
\beta: (U \times U)\times G \to U; & ((x, y), g) \to y,
\\
\epsilon: U \to (U \times U) \times G; & x \to ((x, x), {\frak
e}),
\end{array}
\]
and with multiplication $m: ((U \times U)\times G)_{2} \to (U
\times U)\times G$ and inversion $i: (U \times U) \times G \to (U
\times U) \times G$ defined by
\begin{equation}
\label{mi}
\begin{array}{rcl}
 m(((x, y), g), ((y, z), h))& = &((x, z), g A(x, y)^{-1}h
A(y, z)^{-1} A(x, z)),
\\
i((x, y), g) & = & ((y, x), A(x, y)g^{-1}A(y, x)).
\end{array}
\end{equation}
The fibre over the point $x \in U$ of the Lie algebroid $A((U
\times U)\times G)$ may be identified with the vector space
$T_{x}U \times {\mathfrak g}$. Thus, a section of $A((U \times
U)\times G)$ is a pair $(X, \tilde{\xi})$, where $X$ is a vector
field on $U$ and $\tilde{\xi}$ is a map from $U$ on ${\mathfrak
g}$. Note that the space $\Gamma(A((U \times U)\times G))$ is
generated by sections of the form $(X, 0)$ and $(0, C_{\xi})$,
with $X \in {\mathfrak X}(U)$, $\xi \in {\mathfrak g}$ and
$C_{\xi}: U \to {\mathfrak g}$ being the constant map $C_{\xi}(x)
= \xi$, for all $x \in U$. Moreover, an straightforward
computation, using (\ref{mi}), proves that the vector fields
$\lvec{(X, 0)}$, $\rvec{(X, 0)}$, $\lvec{(0, C_{\xi})}$ and
$\rvec{(0, C_{\xi})}$ on $(U \times U)\times G$  are given by
\begin{equation}
\label{leftright}
\begin{aligned}
 \lvec{(X, 0)}((x, y), g) & =  (0_{x}, X(y), (T_{A(x,
y)}l_{gA(x,y)^{-1}}((T_{y}A_{x})(X(y)))+\\
&\hbox to 3cm{\hfil}{}-(Ad_{A(x,y)^{-1}}(T_yA_y)(X(y)))^l(g))),
\\
\rvec{(X,0)}((x, y), g) & =  (-X(x), 0_y, -(T_{A(x,
y)}l_{gA(x,y)^{-1}}((T_{x}A_{y})(X(x))))+\\
&\hbox to 3cm{\hfil}{}-(Ad_{A(x,y)^{-1}}(T_xA_x)(X(x)))^l(g) )),
\\
\lvec{(0, C_{\xi})}((x, y), g) & =  (0_x, 0_y, (Ad_{A(x,
y)^{-1}}\xi)^l(g)),
\\
\rvec{(0, C_{\xi})}((x, y), g) & =  (0_x, 0_y, \xi^{r}(g)),
\end{aligned}
\end{equation}
for $((x, y), g) \in (U \times U) \times G$, where $l_{h}: G \to
G$ denotes the left-translation in $G$ by $h \in G$, $Ad: G \times
{\mathfrak g} \to {\mathfrak g}$ is the adjoint action of the Lie
group $G$ on ${\mathfrak g}$, $\eta^{l}$ (respectively,
$\eta^{r}$) is the left-invariant (respectively, right-invariant)
vector field on $G$ such that $\eta^{l}({\mathfrak e}) = \eta $
(respectively, $\eta^{r}({\mathfrak e}) = \eta $) and $A_{x}: U
\to G$ and $A_{y}: U \to G$ are the maps defined by
\[
A_{x}(y) = A_{y}(x) = A(x, y).
\]

Now, suppose that $L: (Q \times Q)/G \to \R$ is a Lagrangian
function on the Atiyah groupoid $(Q \times Q)/G$. Then, the
discrete Euler-Lagrange equations for $L$ are
\[
\begin{array}{lcr}
\lvec{(X, 0)}((x, y), g_{k})(L) - \rvec{(X, 0)}((y, z),
g_{k+1})(L) & =& 0,
\\
\lvec{(0, C_{\xi})}((x, y), g_{k})(L) - \rvec{(0, C_{\xi})}((y,
z), g_{k+1})(L) & = & 0,
\end{array}
\]
with $X \in {\mathfrak X}(U)$, $\xi \in {\mathfrak g}$ and $(((x,
y), g_{k}), ((y, z), g_{k+1})) \in ((U \times U)\times G)_{2}$.

From (\ref{leftright}), it follows that the above equations may be
written as
\begin{align}
\label{disLP}
 &D_{2}L ((x, y), g_{k}) + D_{1}L ((y, z), g_{k+1})
 + df_{AL}[x,y, g_{k}](y) +\nonumber\\
&\qquad+df_{AL}[y, z, g_{k+1}](y)
+df^1_{ALI}[x,y,g_k](y)+df^2_{ALI}[y,z,g_k+1](y)= 0,
\\
\label{disLP'} &d(L_{(x, y, \;)}\circ l_{g_{k}} \circ I_{A(x,
y)^{-1}})({\mathfrak e}) - d(L_{(y, z, \;)} \circ
r_{g_{k+1}})({\mathfrak e}) = 0,\qquad\qquad\qquad
\end{align}
where $I_{\bar{g}}: G \to G$ denotes the interior automorphism in
$G$ of $\bar{g}\in G$, $L_{(\bar{x}, \bar{y}, \;)}: G \to \R$ is
the function $L_{(\bar{x}, \bar{y}, \; )}(g) = L(\bar{x}, \bar{y},
g)$, and $f_{AL}[\bar{x}, \bar{y}, \bar{g}]$, $f^1_{ALI}[\bar{x},
\bar{y}, \bar{g}]$ and $f^2_{ALI}[\bar{x}, \bar{y}, \bar{g}]$ are
the real functions on $U$ given by
\begin{align*}
f_{AL}[\bar{x}, \bar{y}, \bar{g}](y)
&= L(\bar{x},\bar{y},\bar{g}A(\bar{x}, \bar{y})^{-1}A(\bar{x}, y)),\\
f^1_{ALI}[\bar{x}, \bar{y}, \bar{g}](y)
&= L(\bar{x},\bar{y},\bar{g}A(\bar{x}, \bar{y})^{-1}A(\bar{y}, y) A(\bar{x},\bar{y})),\\
f^2_{ALI}[\bar{x}, \bar{y}, \bar{g}](y) &=
L(\bar{x},\bar{y},\bar{g}A(\bar{x}, \bar{y})^{-1}A(y,\bar{y})
A(\bar{x},\bar{y})),
\end{align*}
for $\bar{x}, \bar{y}, y \in U $ and $g \in G$. These equations
may be considered as the discrete version of the Lagrange-Poincar\'{e}
equations for a $G$-invariant continuous Lagrangian (see
\cite{CMR} for the local expression of the Lagrange-Poincar\'{e}
equations).

Note that if $A: U \times U \to G$ is the constant map $A(x, y) =
{\mathfrak e}$, for all $(x, y) \in U \times U$, or, in other
words, ${\mathcal A}_{d}$ is the trivial connection then equations
\eqref{disLP} and \eqref{disLP'} may be written as
\begin{equation}
\label{disLP1}
\begin{array}{l}
D_{2}L((x, y), g_{k}) + D_{1}L ((y, z), g_{k+1}) = 0,
\\
\mu_{k+1}(y, z) = Ad^*_{g_{k}} \mu_{k}(x, y),
\end{array}
\end{equation}
where
\[
\mu_{k}(\bar{x}, \bar{y}) = d(r_{g_{k}}^*L_{(\bar{x}, \bar{y},
\;)})({\mathfrak e})
\]
for $(\bar{x}, \bar{y}) \in U \times U$ (compare equations
(\ref{disLP1}) with equations (\ref{5.0*}) and (\ref{5.0'*})).

 %%%%%%%%%%%%%
\subsubsection*{Discrete Elroy's beanie}

As an example of a lagrangian system on an Atiyah grou\-poid, we
consider a discretization of the  Elroy's beanie, which is,
probably,  the most simple example of a dynamical system with a
non-Abelian Lie group of symmetries. The continuous  system
consists in two planar rigid bodies attached at their centers of
mass, moving freely in the plane. The configuration space is
$Q=SE(2)\times S^1$ with coordinates $(x, y, \theta, \psi)$, where
the three first coordinates describe the position and orientation
of the center of mass of the first body and the last one the
relative orientation between both bodies. The continuous system is
described by a Lagrangian $ L_c(x, y, \theta, \psi, \dot{x},
\dot{y}, \dot{\theta}, \dot{\psi})= \frac{1}{2}m
(\dot{x}^2+\dot{y}^2)+\frac{1}{2} I_1\dot{\theta}^2+\frac{1}{2}I_2
(\dot{\theta}+\dot{\psi})^2-V(\psi) $ where $m$ denotes the mass
of the system, $I_1$ and $I_2$ are the inertias of the first and
the second body, respectively, and $V$ is the potential energy.
The system admits reduction by $SE(2)$ symmetry. In fact, the
reduced lagrangian $l_c: TQ/SE(2) \simeq S^{1} \times \R \times
\mathfrak{se}(2) \to \R$ is
\[
l_c(\psi, \dot{\psi}, \Omega_1,
\Omega_2,\Omega_3)=\frac{1}{2}m(\Omega_1^2+\Omega_2^2)+\frac{1}{2}(I_1+I_2)
\Omega_3^2+\frac{1}{2}\frac{I_1I_2}{I_1+I_2}\dot{\psi}^2-V(\psi)
\]
where $\mathfrak{se}(2)$ is the Lie algebra of $SE(2)$,
$\Omega_1=\xi_1$, $\Omega_2=\xi_2$,
$\Omega_3=\xi_3-\frac{I_2}{I_1+I_2}\dot{\psi}$ and  $(\xi_1,
\xi_2, \xi_3)$ are the coordinates of an element of
$\mathfrak{se}(2)$ with respect to the basis $e_1=\left(
\begin{smallmatrix}
0&0&1\\
0&0&0\\
0&0&0
\end{smallmatrix}
\right)$, $e_2=\left(
\begin{smallmatrix}
0&0&0\\
0&0&1\\
0&0&0
\end{smallmatrix}
\right)$ and $e_3=\left(
\begin{smallmatrix}
0&1&0\\
-1&0&0\\
0&0&0
\end{smallmatrix}
\right)$. Note that $\xi_{1} = \dot{x} \cos \theta + \dot{y} \sin
\theta$, $\xi_{2} = -\dot{x} \sin \theta + \dot{y} \cos \theta $
and $\xi_{3} = -\dot{\theta} - \displaystyle \frac{I_{2}}{I_{1} +
I_{2}} \dot{\psi}$ (for more details, see \cite{Lewis,Ostrowski}).

In order to discretize this system, consider $g_k=\left(
\begin{smallmatrix}
\cos\theta_k&-\sin \theta_k&x_k\\
\sin \theta_k&\cos\theta_k&y_k\\
0&0&1 \end{smallmatrix} \right) \in SE(2)$.
 We construct the  discrete connection  ${\mathcal
A}_d:(SE(2)\times S^1)\times (SE(2)\times S^1)\longrightarrow
SE(2)$ defined by ${\mathcal A}_d((g_k, \psi_k), (g_{k+1},
\psi_{k+1}))=g_{k+1}A(\psi_k, \psi_{k+1}) g_{k}^{-1}$, where \[
A(\psi_k, \psi_{k+1})=\left(
\begin{smallmatrix}
\cos(\frac{I_2}{I_1+I_2}\Delta\psi_{k})&-\sin
(\frac{I_2}{I_1+I_2}\Delta\psi_{k})&0\\
\sin (\frac{I_2}{I_1+I_2}\Delta\psi_{k})&\cos(\frac{I_2}{I_1+I_2}\Delta\psi_{k}) &0\\
0&0&1
\end{smallmatrix}
\right)\] Here, $\Delta \psi_{k} = \psi_{k+1} - \psi_{k}$. The
discrete connection ${\mathcal A}_{d}$ precisely induces the
mechanical connection associated with the $SE(2)$-invariant metric
${\mathcal G}$ on $Q$:
\[
{\mathcal G}=m d{x}\otimes dx+ m d{y}\otimes dy+
(I_1+I_2)d\theta\otimes d\theta+I_2 d\theta\otimes d\psi + I_2
d\psi\otimes d \theta+ I_2 d\psi\otimes d\psi
\]

We remark that the continuous Lagrangian $L_{c}$ is the kinetic
energy associated with ${\mathcal G}$ minus the potential energy
$V$.

Next, we consider the Atiyah groupoid ${(Q \times Q)}/{SE(2)}$. As
we know, using the discrete connection ${\mathcal A}_{d}$, one may
define a local isomorphism between the Atiyah groupoid ${(Q \times
Q)}/{SE(2)}$ and the product manifold $U \times U \times SE(2)$,
$U$ being an open subset of $\R$. Then, as a local discretization
of the reduced Lagrangian $l_{c}$, we introduce the discrete
Lagrangian $l_{d}$ on $U \times U \times SE(2)$ given by
\begin{eqnarray*}
&&l_d (\psi_k, \psi_{k+1}, \Omega_{(1)k}, \Omega_{(2)k},
\Omega_{(3)k})= \frac{1}{2h^2}m\left[ \Omega_{(1)k}^2 +
\Omega_{(2)k}^2\right]\\&&\qquad +\frac{(I_1+I_2)}{h^2}
\left[1-\cos(\Omega_{(3)k})\right]
+\frac{1}{2}\frac{I_1I_2}{I_1+I_2}\left(\frac{\Delta
\psi_k}{h}\right)^2-V(\frac{\psi_k+\psi_{k+1}}{2})
\end{eqnarray*}
where $ \Omega_{(1)k}=\Delta x_k\cos\theta_k + \Delta y_{k}\sin
\theta_k,$ $\Omega_{(2)k}=-\Delta x_k\sin\theta_k+\Delta y_k
\cos\theta_k\ $ and $\
\Omega_{(3)k}=-\Delta\theta_{k}-\frac{I_2}{I_1+I_2}\Delta\psi_{k}$.

 Now, if we denote by $\bar{q}_k=(\psi_{k}, \psi_{k+1},
\Omega_{(1)k}, \Omega_{(2)k}, \Omega_{(3)k})$ then
 \begin{eqnarray*}
&\lvec{(0,C_{e_1})}\Big|_{\bar{q}_k}={\small \cos
(\Omega_{(3)k}+\frac{I_2}{I_1+I_2}\Delta\psi_k)\frac{\partial}{\partial
\Omega_{(1)k}} -\sin
(\Omega_{(3)k}+\frac{I_2}{I_1+I_2}\Delta\psi_k)\frac{\partial}{\partial
\Omega_{(2)k}}}&
\\
&\lvec{(0,C_{e_2})}\Big|_{\bar{q}_k}={\small \sin
(\Omega_{(3)k}+\frac{I_2}{I_1+I_2}\Delta\psi_k)\frac{\partial}{\partial
\Omega_{(1)k}} +\cos
(\Omega_{(3)k}+\frac{I_2}{I_1+I_2}\Delta\psi_k)\frac{\partial}{\partial
\Omega_{(2)k}}}&
\\
&\lvec{(0,C_{e_3})}\Big|_{\bar{q}_k}=-\frac{\partial}{\partial
\Omega_{(3)k}},
\quad\rvec{(0,C_{e_1})}\Big|_{\bar{q}_k}=\frac{\partial}{\partial
\Omega_{(1)k}},\quad
\rvec{(0,C_{e_2})}\Big|_{\bar{q}_k}=\frac{\partial}{\partial
\Omega_{(2)k}}&
\\
&\rvec{(0,C_{e_3})}\Big|_{\bar{q}_k}={\small
-\frac{\partial}{\partial
\Omega_{(3)k}}+\Omega_{(2)k}\frac{\partial}{\partial
\Omega_{(1)k}}-\Omega_{(1)k}\frac{\partial}{\partial
\Omega_{(2)k}}} , \quad\lvec{(\frac{\partial}{\partial
\psi},0)}\Big|_{\bar{q}_k}=\frac{\partial }{\partial
\psi_{k+1}},&\\
& \rvec{(\frac{\partial}{\partial
\psi},0)}\Big|_{\bar{q}_k}={\small -\frac{\partial }{\partial
\psi_{k}}+\frac{I_2}{I_1+I_2}\Omega_{(2)k}\frac{\partial}{\partial
\Omega_{(1)k}}-\frac{I_2}{I_1+I_2}\Omega_{(1)k}\frac{\partial}{\partial
\Omega_{(2)k}}}&
\end{eqnarray*}
 Thus, the reduced Discrete Euler-Lagrange equations
\[
\lvec{(0,C_{e_i})}\Big|_{\bar{q}_k}l_d-\rvec{(0,C_{e_i})}\Big|_{\bar{q}_{k+1}}l_d=0,\quad
 \lvec{(\frac{\partial}{\partial
\psi},0)}\Big|_{\bar{q}_k}l_d-\rvec{(\frac{\partial}{\partial
\psi},0)}\Big|_{\bar{q}_{k+1}}l_d=0
\]
are
\[
\left\{
\begin{array}{rcl}
&&\Omega_{(1)k+1}=\Omega_{(1)k} \cos
(\Omega_{(3)k}+\frac{I_2}{I_1+I_2}\Delta\psi_k)- \Omega_{(2)k} \sin (\Omega_{(3)k}+\frac{I_2}{I_1+I_2}\Delta\psi_k)\\
&&\Omega_{(2)k+1}=\Omega_{(1)k} \sin
(\Omega_{(3)k}+\frac{I_2}{I_1+I_2}\Delta\psi_k)+ \Omega_{(2)k} \cos(\Omega_{(3)k}+\frac{I_2}{I_1+I_2}\Delta\psi_k)\\
&&\Omega_{(3)k+1}=\Omega_{(3)k}\\
&&\displaystyle{\frac{I_1I_2}{I_1+I_2}\frac{\psi_{k+2}-2\psi_{k+1}+\psi_{k}}{h^2}}=\displaystyle{-\frac{1}{2}\left(\frac{\partial
V}{\partial \psi}(\frac{\psi_{k+2}+\psi_{k+1}}{2})+\frac{\partial
V}{\partial \psi}(\frac{\psi_{k+1}+\psi_{k}}{2})\right)}
 \end{array}
 \right.
\]
These equations are a discretization of the corresponding reduced
equations for the continuous system (see \cite{Lewis}). In a
forthcoming paper \cite{IMMM}, we will give a complete description
of this example comparing with the continuous equations.
%%%%%%%%%%%%%%%%%%%%%%%%%%%%%%%%%
\subsection{Reduction of discrete Lagrangian
systems}\label{rdLs}

Next, we will present some examples of Lie groupoid epimorphisms
which allow to do reduction.

\vspace{6pt}

\noindent $\bullet$ Let $G$ be a Lie group and consider the pair
groupoid $G\times G$ over $G$. Consider also $G$ as a groupoid
over one point. Then we have that the map
\[
\begin{array}{rccc}
\Phi_l:&G\times G&\longrightarrow&G\\
       &  (g,h)  &\mapsto        &g^{-1}h
\end{array}
\]
is a Lie groupoid morphism, which is obviously a submersion. Thus,
using Corollary \ref{reduction}, it follows that the discrete
Euler-Lagrange equations for a left invariant discrete Lagrangian
on $G\times G$ reduce to the discrete Lie-Poisson equations on $G$
for the reduced Lagrangian. This case appears in~\cite{Mose} as
was first noticed by~\cite{weinstein}, and also appear later
in~\cite{BoSu,BoSu2,MPS1,MPS2}.

Alternatively, one can do reduction of a right-invariant
Lagrangian by using the morphism
\[
\begin{array}{rccc}
\Phi_r:&G\times G&\longrightarrow&G\\
       &  (g,h)  &\mapsto        &gh^{-1}
\end{array}
\]

\noindent $\bullet$ Let $G$ be a Lie group acting on a manifold
$M$ by the left. We consider a discrete Lagrangian on $G\times G$
which depends on the variables of $M$ as parameters $L_m(g,h)$. In
general, the Lagrangian will not be invariant under the action of
$G$, that is $L_m(g,h)\neq L_m(rg,rh)$. Nevertheless, it can
happen that $L_m(rg,rh)=L_{r^{-1}m}(g,h)$. In such cases we can
consider the Lie groupoid $G\times G\times M$ over $G\times M$
where accordingly one consider the elements in $M$ as parameters.
Then the Lagrangian can be considered as a function on the
groupoid $G\times G\times M$ given by $L(g,h,m)\equiv L_m(g,h)$ so
that the above property reads $L(rg,rh,rm)=L(g,h,m)$. Thus we
define the reduction map
\[
\begin{array}{rccc}
\Phi:&G\times G\times M&\longrightarrow&G\times M\\
       &  (g,h,m)  &\mapsto        &(g^{-1}h,g^{-1}m)
\end{array}
\]
where on $G\times M$ we consider the transformation Lie groupoid
defined by the right action $m\cdot g=g^{-1}m$. Since this map is
a submersion, the Euler-Lagrange equations on $G\times G\times M$
reduces to the Euler-Lagrange equations on $G\times M$. This case
occurs in the Lagrange top that was considered as an example in
Section \ref{TaLg} (see also~\cite{BoSu2}).

\vspace{10pt}

\noindent  $\bullet$ Another interesting case is that of a
$G$-invariant Lagrangian $L$ defined on the pair groupoid
$\map{L}{Q\times Q}{\R}$, where $\map{p}{Q}{M}$ is a $G$-principal
bundle. In this case we can reduce to the Atiyah gauge groupoid by
means of the map
\[
\begin{array}{rccc}
\Phi:&Q\times Q&\longrightarrow&(Q\times Q)/G\\
       &  (q,q')  &\mapsto        &[(q,q')]
\end{array}
\]
Thus the discrete Euler-Lagrange equations reduce to the so called
discrete Lagrange-Poincar\'e equations.

\section{Conclusions and outlook}

In this paper we have elucidated the geometrical framework for
discrete Mechanics on Lie groupoids. Using as a main tool the
natural Lie algebroid structure on the vector bundle $\pi^\tau:
{\mathcal P}^{\tau}G \to G$ we have found intrinsic expressions
for the discrete Euler-Lagrange equations. We introduce the
Poincar\'e-Cartan sections,  the discrete Legendre transformations
and the discrete evolution operator in the Lagrangian and in the
Hamiltonian formalism. The notion of regularity has been
completely characterized and we prove the symplecticity of the
discrete evolution operators. Moreover,  we have studied the
symmetries of discrete Lagrangians on Lie groupoids relating them
with constants of the motion via Noether's Theorems. The
applicablity of these developments has been stated in several
interesting examples, in particular for the case of discrete
Lagrange-Poincar\'e equations. In fact, the general theory of
discrete symmetry reduction naturally follows from our results.

In this paper we have confined ourselves to the geometrical
aspects of mechanics on Lie groupoids. In a forthcoming paper (see
\cite{IMMM}) we will study the construction of geometric
integrators for mechanical systems on Lie algebroids. We will
introduce the exact discrete Lagrangian and we will discuss
different discretizations of a continuous Lagrangian and its
numerical implementation.

Another different aspect we will work on it in the future is to
develop natural extensions of the above theories for forced
systems and systems with holonomic and nonholonomic constraints.

\end{document}